\documentclass{amsart}
\usepackage{amssymb,mathrsfs,amsmath,dsfont,array}
\usepackage{multirow,arydshln}
\usepackage{tikz}
\usetikzlibrary{arrows,positioning}
\usepackage{capt-of}
\usepackage{mathbbol}     
\usepackage{color}
\usepackage{verbatim} 
\usepackage[all,cmtip]{xy}
\usepackage{amsbsy}
\usepackage{amsmath}
\usepackage{nccmath}

\usepackage{hyperref}

\newtheorem{Thm}{Theorem}[section]
\newtheorem{Pro}[Thm]{Proposition}
\newtheorem{cor}[Thm]{Corollary}
\newtheorem{lem}[Thm]{Lemma}

\newtheorem{deft}[Thm]{Definition}
\newtheorem{rem}[Thm]{Remark}

\makeatletter
\newcommand{\AlignFootnote}[1]{%
    \ifmeasuring@
    \else
        \footnote{#1}%
    \fi
}
\makeatother
\numberwithin{equation}{section}

\newcommand{\End}{\mathop{{\rm End}}\nolimits}

\renewcommand{\tilde}{\widetilde}

\def\a1s{a_1,\cdots, a_s}
\def\a{\alpha}

\def\aa{\mathcal A}
\def\aA{\mathscr{A}}

\def\andd{\quad\hbox{and}\quad}

\def\ad{\hbox{ad}}

\def\b{\beta}

\def\bl4{B_{\ell\geq4}}

\def\bbbc{{\mathbb C}}

\def\d{\delta}

\def\gG{\mathscr{G}}
\def\gG{{\mathcal G}}

\def\fg{\mathfrak{g}}

\def\scg{\mathscr{G}}

\def\hH{\mathscr{H}}

\def\hh{{\mathcal H}}

\def\fh{\mathfrak{h}}
\def\sch{\mathscr{H}}

\def\fk{\mathfrak{k}}

\def\lam{\lambda}
\def\Lam{\Lambda}
\def\LL{\mathcal{L}}
\def\scl{\mathscr{L}}

\def\fl{\mathfrak{L}}

\def\ep{\epsilon}
\def\fm{(\cdot,\cdot)}

\def\bbbr{{\mathbb R}}

\def\supp{\hbox{\rm supp}}

\def\1k{\frac{1}{k}}
\def\op{\oplus}
\def\ot{\otimes}

\def\sub{\subseteq}
\def\sg{\sigma}

\def\pf{\noindent{\bf Proof. }}

\def\sspan{\hbox{\rm span}}

\def\i{{\mathcal I}}

\def\ft{\mathfrak{t}}

\def\u{{\mathcal U}}

\def\v{{\mathcal V}}

\def\w{{\mathcal W}}

\def\bbbz{{\mathbb Z}}

\def\1il{1\leq i\leq\ell}

\begin{document}
\title{Finite weight modules over twisted affine Lie superalgebras}
\thanks{2010 Mathematics Subject Classification: 17B10, 17B65.}
\thanks{Key Words: Finite Weight Modules; Twisted Affine Lie Superalgebras.}
\maketitle
\centerline{Malihe Yousofzadeh\footnote{Department of Pure Mathematics, Faculty of Mathematics and Statistics, University of Isfahan, Isfahan, Iran,
P.O.Box 81746--73441, and School of Mathematics, Institute for Research in
Fundamental
Sciences (IPM), P.O. Box: 19395-5746, Tehran, Iran.
Email address: ma.yousofzadeh@sci.ui.ac.ir \& ma.yousofzadeh@ipm.ir.\\This research  was in part supported by
a grant from IPM (No. 98170424) and  is partially carried out in
IPM-Isfahan Branch.}}

\begin{abstract}
This work provides  the first step toward the classification of irreducible finite weight modules over twisted affine Lie superalgebras. We divide the class of such modules into two subclasses called hybrid and tight. We reduce the classification of  hybrid irreducible finite weight modules to the classification of  cuspidal modules of finite dimensional  cuspidal {Lie superalgebras} which is discussed in a work of Dimitrov, Mathieu and Penkov.
\end{abstract}

\section{Introduction}
To state the results of this paper, we need to start with some definitions. Suppose that $\scl=\scl_0\op\scl_1$ is a
Lie superalgebra with a splitting Cartan subalgebra $\sch\sub\scl_0$ and corresponding root system $R.$
An $\mathscr{L}$-module $M$ is said to have a weight space decomposition with respect to $\mathscr{H}$ (or a weight module) if $$M=\op_{\lam\in \mathscr{H}^\ast}M^\lam$$ in which $\mathscr{H}^\ast$ is the dual space of $\mathscr{H}$  and
$$M^\lam:=\{v\in M\mid hv=\lam(h)v\;\;(h\in\mathscr{H})\}\quad\quad(\lam\in \mathscr{H}^\ast).$$ If each  $M^\lam$ is finite dimensional, the module $M$ is called a finite weight module.
To study the weight modules over $\mathscr{L},$ some subsets of $R$  satisfying  $(P+P)\cap R\sub P$ and $R=P\cup -P,$ get involved; such subsets are called parabolic subsets.
 For a parabolic subset $P$ of $R,$ we  have the decomposition
\[\mathscr{L}=\mathscr{L}^+\op\mathscr{L}^\circ\op\mathscr{L}^-\] where
 $$\hbox{\small $\mathscr{L}^\circ:=\op_{\a\in P\cap -P}\mathscr{L}^\a,\; \mathscr{L}^+:=\op_{\a\in P\setminus-P}\mathscr{L}^\a\andd \mathscr{L}^-:=\op_{\a\in - P\setminus P}\mathscr{L}^\a.$}$$
We  set $$\frak{p}:=\mathscr{L}^\circ\op\mathscr{L}^+.$$

 For a  functional $\lam$ on the $\bbbr$-linear span of $R,$ we have the decomposition  $R=R^+\cup R^\circ\cup R^-,$ called a triangular decomposition, where $$R^\pm:=\{\a\in R\mid \lam(\a)\gtrless0\} \andd R^\circ:=\{\a\in R\mid \lam(\a)=0\}.$$ In this case, $P_\lam:=R^+\cup R^\circ$ is a parabolic subset of $R.$ Moreover, if $\mu $ is a functional  on the $\bbbr$-linear span of $R^\circ,$ we have a triangular decomposition $R^\circ=R^{\circ,+}\cup R^{\circ,\circ}\cup R^{\circ,-}$ for $R^\circ$ and  $P_{\lam,\mu}:=R^+\cup R^{\circ,+}\cup R^{\circ,\circ}$ is also a parabolic subset of $R.$ We note that $P_{\lam,0}=P_{\lam}.$

For functionals $\lam$ and $\mu$ as above, consider subalgebras     $\scl^\circ$ and $\frak{p}$ corresponding to $P_{\lam,\mu}.$ Each irreducible $\scl^\circ$-module $N$ is a module of $\frak{p}$ with trivial action of $\mathscr{L}^+.$ Then  $$\tilde N:=U(\mathscr{L})\ot_{U(\mathfrak{p})}N$$ is an $\mathscr{L}$-module; here  $U(\mathscr{L})$ and $U(\frak{p})$ denote respectively  the universal enveloping algebras of $\mathscr{L}$ and $\frak{p}.$ If the $\scl$-module  $\tilde N$ contains a  maximal submodule $Z$ {intersecting} $N$ trivially, the quotient module $${\rm Ind}_{\mathscr{L}}(N):=\tilde N/Z$$ is called  a parabolically induced module if $\lam$ is nonzero.
An irreducible $\scl$-module  which is not parabolically induced is called cuspidal.

The study of finite weight modules of Lie (super)algebras has an ancient root in the literature. In  \cite{BL1}, \cite{BL2} and \cite{F}, the authors   classify  irreducible finite weight modules of finite dimensional reductive Lie algebras. The important point to get this classification is that the classification is reduced to the classification of cuspidal modules.

This perspective can be developed to current Lie (super)algebras, finite dimensional basic classical simple Lie superalgebras and affine Lie (super)algebras; see \S \ref{review} for the  review of the literature.

Suppose $\mathscr{L}$ is a twisted affine Lie superalgebra of type $X=A(2k-1,2\ell-1)^{(2)}$ ({\tiny$(k,\ell)\neq (1,1)$}), $A(2k,2\ell)^{(4)},$    $A(2k,2\ell-1)^{(2)}$ and $D(k+1,\ell)^{(2)}$ where  $k,\ell$ are positive integers, with standard Cartan subalgebra $\hH.$
 The  root system $R$ of $\mathscr{L}$ with respect to $\mathscr{H}$ has three kind of roots: nonzero real roots (roots which are not self-orthogonal with respect to the canonical bilinear form on the dual space of $\hH$), imaginary  roots (roots which are orthogonal to all roots) and nonsingular roots (neither real nor imaginary). Nonsingular roots appear just  as the weights for the $\mathscr{H}$-module $\mathscr{L}_1$ and all roots of the $\mathscr{H}$-module $\mathscr{L}_0$ are real but the odd part $\scl_1$ may contain real roots as well.
Due to the existence of roots which are either nonsingular or odd real,  representation theory of affine Lie superalgebras is more complicated  comparing with the non-super case.

We next suppose  $M$ is an  irreducible finite weight  module over the twisted affine Lie superalgebra $\mathscr{L}$. Then, each nonzero root vector corresponding to a nonzero real root $\a$, acts on $M$ either injectively or locally nilpotently.
We denote by $R^{in}$ (resp. $R^{ln}$), the subset of $R$ consisting of all nonzero real roots whose nonzero root vectors act injectively (resp. locally nilpotently).
 If $R^{ln}$ coincides with the set $R_{re}^\times$ of all nonzero real roots, then   $M$  is called integrable.

We know that the imaginary roots of the twisted affine Lie superalgebra $\mathscr{L}$ generates a free abelian group of rank 1; say e.g., $\bbbz\d.$ We  show that for each nonzero real root $\a,$ one of the following occurs:
\begin{itemize}
\item  $\a$ is full-locally nilpotent, i.e., $R\cap (\a+\bbbz\d)\sub R^{ln},$
\item $\a$ is full-injective, i.e.,  $R\cap (\a+\bbbz\d)\sub R^{in},$
\item  $\pm\a$ are up-nilpotent hybrid, i.e., there is a positive integer $m$ with $$R\cap (\pm\a+\bbbz^{\geq m}\d)\sub R^{ln} \andd R\cap (\pm\a+\bbbz^{\leq -m}\d)\sub R^{in},$$
\item $\pm\a$ are down-nilpotent hybrid, i.e., there is a positive integer $m$ with $$R\cap (\pm\a+\bbbz^{\geq m}\d)\sub R^{in}\andd R\cap (\pm\a+\bbbz^{\leq -m}\d)\sub R^{ln}.$$
\end{itemize}
{Up to a weight $\hH$-module whose weights are  nonzero imaginary roots}, the even part of  $\mathscr{L}$ is  a summation of two affine Lie algebra  $\scg_1$ and $\scg_2$ with  corresponding root systems $R(1)$ and $R(2)$ respectively.
{We call the irreducible finite weight   $\mathscr{L}$-module $M$  hybrid if all nonzero real roots of $R(1)$ and $R(2)$ are hybrid and otherwise call it tight.
}
If $i\in\{1,2\}$ and all nonzero real roots of $R(i)$  are  hybrid, then either all of them are up-nilpotent hybrid  or all of them are down-nilpotent hybrid. We show that there exists a compatibility between $R(1)$ and $R(2);$ i.e., we prove that  if  all nonzero real roots of $R(1)\cup R(2)$ are hybrid, then either all of them are up-nilpotent hybrid  or all of them are down-nilpotent hybrid. Having this in hand, we then get a nontrivial triangular decomposition  $R^+\cup R^\circ\cup R^-$ for $R$ in case $M$ is hybrid. The next step is finding nonzero weight vectors $v$ with $\mathscr{L}^\a v=\{0\}$ for all  $\a\in R^+.$ Since $R_{re}^\times=R^{ln}\cup R^{in},$  we can show that there are nonzero weight vectors $v$ with $\mathscr{L}^\a v=\{0\}$ for all  real roots $\a\in R^+$ whether odd or even and also for all imaginary roots $\a\in R^+.$
We then go through the nonsingular roots of $R^+;$ more precisely, among   nonzero weight vectors $v$ with $\mathscr{L}^\a v=\{0\}$ for all  real and imaginary roots $\a\in R^+,$ we find some satisfying $\mathscr{L}^\a v=\{0\}$ for all  nonsingular roots $\a\in R^+.$ This shows that $$M^{\mathscr{L}^+}:=\{v\in M\mid \mathscr{L}^\a v=\{0\}\;\;(\a\in R^+)\}$$ is a nonzero irreducible finite weight $\mathscr{L}^\circ$-module and $M$ is parabolically induced from $M^{\mathscr{L}^+}.$ Moreover,  we prove that if $M$ is hybrid,
the classification problem is reduced to
the classification of irreducible finite weight
cuspidal modules over finite-dimensional cuspidal Levi subsuperalgebras discussed
by Dimitrov, Mathieu and Penkov \cite{DMP}.

The outline of the paper is as follows:  After ``Introduction'' and ``Review of The Literature'', in Section 3, we first gather some information regarding twisted affine Lie superalgebras of types $X=A(2k-1,2\ell-1)^{(2)}$ ({\tiny$(k,\ell)\neq (1,1)$}), $A(2k,2\ell)^{(4)},$    $A(2k,2\ell-1)^{(2)}$ and $D(k+1,\ell)^{(2)}$ where  $k,\ell$ are positive integers and then  prove   general information regarding weight modules. In Section 4, we focus on modules having shadow; see Definition~\ref{shadow}. Section 5 is devoted to our main results. We end up the paper with an appendix section in which, for the convenience of readers, we recall the structure of twisted affine Lie superalgebras.
\section{Review of the literature}\label{review}
In this section, we give a history of the study of finite weight modules of Lie (super)algebras.
Suppose that $R$ is the root system of a Lie superalgebra $\scl=\scl_0\op\scl_1$ with respect to a splitting Cartan subalgebra $\sch\sub \scl_0$ and $M$ is an irreducible finite weight $\scl$-module.

If $\mathscr{L}$ is a finite dimensional  reductive Lie algebra and     both  $R^{ln}$ and $R^{in}$ are nonempty subsets of $R^\times$, then  $P:=R^{ln}\cup -R^{in}\cup \{0\}$ is a parabolic subset of $R$. This in turn implies that there is a functional $\lam$ on the $\bbbr$-linear span of $R$ such that $P=R^+\cup R^\circ$ \cite[Pro. VI.7.20]{Bo}. Then it follows  that
$M^{\mathscr{L}^+}$ is an irreducible finite weight  $\mathscr{L}^\circ$-module and $M$ is isomorphic to the module which is  parabolically induced  from $M^{\mathscr{L}^+}.$ The $\mathscr{L}^\circ$-module $M^{\mathscr{L}^+}$ is a tensor product of a {finite dimensional module} and a  finite weight module on which all nonzero roots act injectively; in fact a cuspidal module; see \cite[Thm.~4.18]{F} and \cite[Cor.~3.7]{DMP}.

In   affine Lie algebra case, the existence of imaginary roots (i.e., those roots which are orthogonal to all  roots) makes the study  more complicated.
An affine Lie algebra $\mathscr{L}$ has a 1-dimensional center $\bbbc c.$ The central element  $c$ acts on the irreducible $\mathscr{L}$-module $M$ as $\lam {\rm id}.$ This $\lam$ is called the level of $M.$    In \cite{C}--\cite{CP2}, the authors study  integrable irreducible finite weight modules over affine Lie algebras; to study zero level modules, they introduce certain modules called loop modules. Irreducible finite weight loop modules are classified  in \cite{E}.  Then in \cite{F1}--\cite{F5} and \cite{FS}, the authors   study nonzero level  irreducible finite weight modules over affine Lie algebras.

 Each affine root system is a subset of $\dot R+\bbbz\d$ where $\dot R$ is an irreducible finite root system and  $\d$ is an imaginary root such that  $\bbbz\d$ is  the group generated by the imaginary roots.  The following two cases can happen:
\begin{itemize}
\item  for all  $\dot\a\in \dot R^\times,$  both sets  $R^{ln}\cap (\dot\a+\bbbz\d)$ and  $R^{in}\cap (\dot\a+\bbbz\d)$ are nonempty,
\item there exists  $\dot\a\in \dot R^\times$ such that   $R\cap (\dot\a+\bbbz\d)\sub R^{ln}$ or  $R\cap (\dot\a+\bbbz\d)\sub R^{in}.$
\end{itemize}
The authors in \cite{DG} show that in the former case, either $P:=R^{ln}\cup -R^{in}\cup \bbbz^{\geq 0}\d$ or $P:=R^{ln}\cup -R^{in}\cup \bbbz^{\leq 0}\d$ is a parabolic subset of $R$ and in the latter case
for
\begin{align*}
\hbox{\small$\dot R^{i}:=\{\dot\a\in \dot R_{re}\mid (\dot\a+\bbbz\d)\cap R\sub R^{in} \},\;
\dot R^{f}:=\{\dot\a\in \dot R_{re}\mid (\dot\a+\bbbz\d)\cap R\sub R^{ln} \}$}
\end{align*} and
$ \hbox{\small$\dot R^m:=\dot R\setminus(\dot R^{i}\cup \dot R^{f})$},$ the set
$$P:=((\dot R^{f}\cup -\dot R^{i}\cup \dot R^{m})+\bbbz\d)\cap R$$ is a parabolic subset of $R.$ Using the identification of parabolic subsets in \cite{DFG}, $P=R^+\cup R^\circ$ for a triangular decomposition $R=R^+\cup R^\circ\cup R^-.$ This helps them  to prove that if $R^{ln}$ is a nonempty proper subset of the set of  nonzero real roots $R_{re}^\times,$ then $M^{\mathscr{L}^+}$ is  an irreducible module of $\mathscr{L}^\circ$   and that $M$ is isomorphic to the module which is parabolically induced  from $M^{\mathscr{L}^+}.$   Then they study those irreducible  finite weight modules with $R^{in}=R_{re}^\times$.

In 2001, I. Dimitrov and his coauthors initiated the study of infinite dimensional irreducible  finite weight modules of Lie superalgebras \cite{DMP}. They classified irreducible finite weight modules of basic classical simple Lie superalgebras by reducing the classification problem to the classification of cuspidal modules. Then in 2006,  S. Eswara Rao and V. Futorny \cite{F6}, \cite{EF} classified {irreducible} finite weight modules over untwisted affine Lie superalgebras on which the canonical central element acts as a nonzero multiple of the identity map. Recently, L. Calixto and  V. Futorny have studied highest weight modules over untwisted  affine Lie superalgebras \cite{CF}.
In this work, we  continue the  study of finite weight modules; we study finite weight modules   over twisted affine Lie superalgebras $A(2k-1,2\ell-1)^{(2)}$ ({\tiny$(k,\ell)\neq (1,1)$}), $A(2k,2\ell)^{(4)},$    $A(2k,2\ell-1)^{(2)}$ and $D(k+1,\ell)^{(2)}$ where  $k,\ell$ are positive integers. We complete the study of hybrid modules and pave the way to start the study of tight irreducible finite weight modules. In an ongoing paper, we are dealing with  irreducible (weak) integrable finite weight  modules.

\section{Generic weight modules}\label{generic}
Throughout this section, we assume  $\fl=\fl_0\op\fl_1$ is a {twisted affine Lie superalgebra  of type $X=A(2k-1,2\ell-1)^{(2)}$ ({\tiny$(k,\ell)\neq (1,1)$}), $A(2k,2\ell)^{(4)},$    $A(2k,2\ell-1)^{(2)}$ and $D(k+1,\ell)^{(2)}$ in which $k,\ell$ are positive integers; {see Appendix for the details regarding the structure of twisted affine Lie superalgebras}. Suppose that $\fh\sub \fl_0$ is the standard Cartan subalgebra of $\fl$ with corresponding root system $R.$}
We mention that   $R=R_0\cup R_1$ where $R_0$ (resp. $R_1$)   is the set of weights of $\fl_0$ (resp. $\fl_1$) with respect to $\fh.$

  One also knows that $\fl$ is equipped with a nondegenerate  (super)symmetric invariant bilinear form $\fm.$ As the form is nondegenerate on  $\fh,$ one can transfer the  form  on  $\fh$ to a form on $\fh^*$ denoted again by $\fm.$ We set
\begin{equation}\label{im-ns}
\begin{array}{lll}
R_{re}^\times := \{\a\in R\mid (\a,\a)\neq 0\},& R_{re}:=\{0\}\cup R_{re}^\times & \hbox{(real roots)},\\
R_{im}:=\{\a\in R\mid (\a,\b)=0\;\;\forall \b\in R\},& R_{im}^\times:=R_{im}\setminus\{0\}& \hbox{(imaginary  roots)},\\
R_{ns}:=\{0\}\cup (R\setminus (R_{re}\cup R_{im})),& R_{ns}^\times:=R_{ns}\setminus\{0\}& \hbox{(nonsingular  roots)}.
\end{array}
\end{equation}
It is known that  $R_{im}$ generates a free abelian group of rank 1; say $\bbbz\d.$ Also,
\begin{equation}
\label{dim}
\dim(\fl^\a)=1\quad\quad(\a\in R\setminus R_{im})
\end{equation}
and
\begin{equation}
\label{sl2}
\parbox{2.8in}{if $\a\in R_{re}^\times\cap R_0,$ then there are $e\in\fl^\a$ and $f\in\fl^{-\a}$ such that $(e,f,[e,f])$ is an $\frak{sl}_2$-triple.}
\end{equation}
The root system $R$   has an expression  as   in the following table:

  \begin{table}[h]\caption{\small Root systems of twisted affine Lie superalgebras} \label{table1}
 {\footnotesize \begin{tabular}{|c|l|}
\hline
$X^{({m})}$ &\hspace{3.25cm}$R$ \\
\hline
$A(2k,2\ell-1)^{(2)}$&$\begin{array}{rcl}
\bbbz\d
&\cup& \bbbz\d\pm\{\ep_i,\d_j,\ep_i\pm\ep_r,\d_j\pm\d_s,\ep_i\pm\d_j\mid i\neq r,j\neq s\}\\
&\cup& (2\bbbz+1)\d\pm\{2\ep_i\mid 1\leq i\leq k\}\\
&\cup& 2\bbbz\d\pm\{2\d_j\mid 1\leq j\leq \ell\}.
\end{array}$\\
\hline
$A(2k-1,2\ell-1)^{(2)},\;(k,\ell)\neq (1,1)$& $\begin{array}{rcl}
\bbbz\d&\cup& \bbbz\d\pm\{\ep_i\pm\ep_r,\d_j\pm\d_s,\ep_i\pm\d_j\mid i\neq r,j\neq s\}\\
&\cup& (2\bbbz+1)\d\pm\{2\ep_i\mid 1\leq i\leq k\}\\
&\cup& 2\bbbz\d\pm\{2\d_j\mid 1\leq j\leq \ell\}
\end{array}$\\
\hline
$A(2k,2\ell)^{(4)}$& $\begin{array}{rcl}
\bbbz\d&\cup&  \bbbz\d\pm\{\ep_i,\d_j\mid 1\leq i\leq k,\;1\leq j\leq \ell\}\\
&\cup& 2\bbbz\d\pm\{\ep_i\pm\ep_r,\d_j\pm\d_s,\d_j\pm\ep_i\mid i\neq r,j\neq s\}\\
&\cup&(4\bbbz+2)\d\pm\{2\ep_i\mid 1\leq i\leq k\}\\
&\cup& 4\bbbz\d\pm\{2\d_j\mid 1\leq j\leq \ell\}
\end{array}$\\
\hline
$D(k+1,\ell)^{(2)}$& $\begin{array}{rcl}
\bbbz\d&\cup&  \bbbz\d\pm\{\ep_i,\d_j\mid 1\leq i\leq k,\;1\leq j\leq \ell\}\\
&\cup& 2\bbbz\d\pm\{2\d_j,\ep_i\pm\ep_r,\d_j\pm\d_s,\d_j\pm\ep_i\mid i\neq r,j\neq s\}
\end{array}$\\
\hline
 \end{tabular}}
 \end{table}
with
\[R_{ns}^\times=R\cap (\bbbz\d\pm\{\ep_i\pm\d_j\mid 1\leq i\leq k,\;1\leq j\leq \ell\}).\] One can see that
\begin{equation}\label{re-re}
(R^\times_{ns}+R^\times_{ns})\cap R\sub R_{re}\cup R_{im}.
\end{equation}

\newpage

The root system $R_0$ of $\LL_0$ is as follows:

\begin{table}[h]\caption{\small The zero part of the root systems}
{\footnotesize \begin{tabular}{|c|l|}
\hline
$X^{(m)}$ &\hspace{3cm}$R_{0}$\\
\hline
$A(2k,2\ell-1)^{(2)}$&$\begin{array}{rcl}
\bbbz\d
&\cup& \bbbz\d\pm\{\ep_i,\ep_i\pm\ep_r,\d_j\pm\d_s\mid i\neq r,j\neq s\}\\
&\cup& (2\bbbz+1)\d\pm\{2\ep_i\mid 1\leq i\leq k\}\\
&\cup& 2\bbbz\d\pm\{2\d_j\mid 1\leq j\leq \ell\}.
\end{array}$
\\
\hline
$\begin{array}{c}
A(2k-1,2\ell-1)^{(2)}\\(k,\ell)\neq (1,1)
\end{array}
$& $\begin{array}{rcl}
\bbbz\d&\cup& \bbbz\d\pm\{\ep_i\pm\ep_r,\d_j\pm\d_s\mid i\neq r,j\neq s\}\\
&\cup& (2\bbbz+1)\d\pm\{2\ep_i\mid 1\leq i\leq k\}\\
&\cup& 2\bbbz\d\pm\{2\d_j\mid 1\leq j\leq \ell\}
\end{array}$
\\
\hline
$A(2k,2\ell)^{(4)}$& $\begin{array}{rcl}
2\bbbz\d&\cup&  2\bbbz\d\pm\{\ep_i\mid 1\leq i\leq k\}\\
&\cup& (2\bbbz+1)\d\pm\{\d_j\mid 1\leq j\leq \ell\}\\
&\cup& 2\bbbz\d\pm\{\ep_i\pm\ep_r,\d_j\pm\d_s\mid i\neq r,j\neq s\}\\
&\cup&(4\bbbz+2)\d\pm\{2\ep_i\mid 1\leq i\leq k\}\\
&\cup& 4\bbbz\d\pm\{2\d_j\mid 1\leq j\leq \ell\}
\end{array}$
\\
\hline
$
D(k+1,\ell)^{(2)}
$& $\begin{array}{rcl}
\bbbz\d&\cup&  \bbbz\d\pm\{\ep_i\mid 1\leq i\leq k,\;1\leq j\leq \ell\}\\
&\cup& 2\bbbz\d\pm\{2\d_j,\ep_i\pm\ep_r,\d_j\pm\d_s\mid i\neq r,j\neq s\}
\end{array}$
\\
\hline
 \end{tabular}}
 \end{table}

\noindent We  see that
\begin{equation}
\label{span}
\sspan_\bbbr R_0=\sspan_\bbbr R=\sspan_\bbbr\{\d,\ep_i,\d_j\mid 1\leq i\leq k,1\leq j\leq \ell\}.
\end{equation}
Also,
 there is a positive integer $r$ with
 \begin{align}\label{Salpha}
R_0+r\bbbz\d\sub R_0\andd R_1+r\bbbz\d\sub R_1.
\end{align}

We also have from Table~\ref{table1} that  $R\sub \dot R+\bbbz\d$ where $\dot R$ is as in the following table:

{\tiny
\begin{table}[h]\caption{\small $R$ modulo $\bbbz\d$}\label{dot}
{\footnotesize\begin{tabular}{|c|c|}
\hline
$X^{(m)}$ & $\dot R$ \\
\hline
$A(2k,2\ell-1)^{(2)}$& $\pm\{\ep_i,\d_j,\ep_i\pm\ep_r,\d_j\pm\d_s,\ep_i\pm\d_j\mid 1\leq i,r\leq k,\; 1\leq j,s\leq \ell\}$\\
\hline
$A(2k-1,2\ell-1)^{(2)}$& $\pm\{\ep_i\pm\ep_r,\d_j\pm\d_s,\ep_i\pm\d_j\mid 1\leq i,r\leq k,\; 1\leq j,s\leq \ell\}$\\
$(k,\ell)\neq (1,1)$&\\
\hline
$A(2k,2\ell)^{(4)}$& $\pm\{\ep_i,\d_j,\ep_i\pm\ep_r,\d_j\pm\d_s,\ep_i\pm\d_j\mid 1\leq i,r\leq k,\; 1\leq j,s\leq \ell\}$\\
\hline
$D(k+1,\ell)^{(2)}$&$\pm\{\ep_i,\d_j,\ep_i\pm\ep_r,\d_j\pm\d_s,\ep_i\pm\d_j\mid 1\leq i\neq r\leq k,\; 1\leq j,s\leq \ell\}$\\
\hline
 \end{tabular}
}
 \end{table}
 }

An element $\dot\a\in\dot R$ is called real (resp. nonsingular) if it is either $0$ or  $(\dot\a+\bbbz\d)\cap R\sub R_{re}$ (resp. $R_{ns}$).
The set $\dot R_{re}$ of real roots of $\dot R$ is a finite root system with a decomposition  $\dot R_{re}=\dot R_1\cup \dot R_2$ into two irreducible finite root systems
$\dot R_1$ and $\dot R_2$. We set
\begin{equation}\label{new11}
\dot R_*=(\dot R_1)_*\cup (\dot R_2)_*\quad\quad\quad (*=sh,lg,ex);
\end{equation}
 here ``~$sh$~'', ``~$lg$~'' and ``~$ex$~'' stand respectively for short, long and extra long roots.
Setting
\begin{equation}\label{s}T_{\dot\a}:=\{\sg\in\bbbz\d\mid \dot\a+\sg\in R_0\}\andd S_{\dot\a}:=\{\sg\in\bbbz\d\mid \dot\a+\sg\in R\}\quad\quad(\dot\a\in \dot R),\end{equation} we get

\newpage
\begin{table}[h]\caption{\small Extensions of the elements of $\dot R$} \label{table2}
\footnotesize{\begin{tabular}{|l|c|c|c|c|}
  \hline
  & $A(2k,2\ell-1)^{(2)}$& $A(2k-1,2\ell-1)^{(2)}$&$A(2k,2\ell)^{(4)}$&$D(k+1,\ell)^{(2)}$\\
  \hline
  $S_{\pm\ep_i}$&$\bbbz\d$&{$\emptyset$}&$\bbbz\d$&$\bbbz\d$\\
  \hline
  $S_{\pm\ep_i\pm\ep_j}$&$\bbbz\d$&$\bbbz\d$&$2\bbbz\d$&$2\bbbz\d$\\
  \hline
  $S_{\pm2\ep_i}$&$(2\bbbz+1)\d$&$(2\bbbz+1)\d$&$(4\bbbz+2)\d$&{$\emptyset$}\\
  \hline
  $S_{\pm\d_j}$&$\bbbz\d$&{$\emptyset$}&$\bbbz\d$&$\bbbz\d$\\
  \hline
  $S_{\pm\d_j\pm\d_s}$&$\bbbz\d$&$\bbbz\d$&$2\bbbz\d$&$2\bbbz\d$\\
  \hline
  $S_{\pm2\d_j}$&$2\bbbz\d$&$2\bbbz\d$&$4\bbbz\d$&$2\bbbz\d$\\
  \hline
  $S_{\pm\ep_i\pm\d_j}$&$\bbbz\d$&$\bbbz\d$&$2\bbbz\d$&$2\bbbz\d$\\
  \hline
  \hline
  $T_{\pm\ep_i}$&$\bbbz\d$&{$\emptyset$}&$2\bbbz\d$&$\bbbz\d$\\
  \hline
  $T_{\pm\ep_i\pm\ep_j}$&$\bbbz\d$&$\bbbz\d$&$2\bbbz\d$&$2\bbbz\d$\\
  \hline
  $T_{\pm2\ep_i}$&$(2\bbbz+1)\d$&$(2\bbbz+1)\d$&$(4\bbbz+2)\d$&{$\emptyset$}\\
  \hline
  $T_{\pm\d_j}$&{$\emptyset$}&{$\emptyset$}&$(2\bbbz+1)\d$&{$\emptyset$}\\
  \hline
  $T_{\pm\d_j\pm\d_s}$&$\bbbz\d$&$\bbbz\d$&$2\bbbz\d$&$2\bbbz\d$\\
  \hline
  $T_{\pm2\d_j}$&$2\bbbz\d$&$2\bbbz\d$&$4\bbbz\d$&$2\bbbz\d$\\
  \hline
\end{tabular}}
\end{table}
 One can easily see  from this table  that setting  $\mathcal{R} $ to be  either $R_0$ or $R$ and $\dot{\mathcal{R}}$ to be respectively  $\dot R_0:=\{\dot\gamma\in \dot R\mid (\dot\gamma+\bbbz\d)\cap R_0\neq \emptyset\}$ or $\dot R,$ then
 \begin{equation}
 \label{imp}
 \parbox{4.2in}{for $0\neq \dot\a\in \dot{\mathcal{R}},$
 $\{m\d\mid \dot\a+m\d\in \mathcal{R}\}=(r_{\dot\a}\bbbz+k_{\dot\a})\d$ for some
 $r_{\dot\a}\in\{1,2,4\}$
 and $0\leq k_{\dot\a}<r_{\dot\a}.$  Moreover,  there is
 $0\neq \dot \a^*\in \dot{\mathcal{R}}_{re}$ with $\{m\d\mid \dot\a^*+m\d\in \mathcal{R}\}=r_{\dot\a^*}\bbbz\d=\mathcal{R}_{im}$ and $r_{\dot\a^*}\mid r_{\dot\a},$ for all $\dot\a\in \dot R^\times.$
 }
 \end{equation}

\begin{rem}\label{decom-aff}
{\rm
If  $\fl=A(2k-1,2\ell-1)^{(2)},$ $(k,\ell)\neq (1,1)$, then $R_{re}\sub R_0,$ so
\begin{equation}\label{new10}
R\cap (R_{ns}^\times+R_{re}^\times)\sub R_{ns}^\times
\end{equation}
 as $[\fl_1^\a,\fl_0^\b]\sub \fl_1^{\a+\b}$
for  $\a\in R_{ns}^\times$ and $\b\in R_{re}^\times.$
\smallskip

Also as $(k,\ell)\neq(1,1),$ without loss of generality, we assume $\ell>1.$
For $\dot\a_1,\dot\a_2,\dot\a_3\in\{\pm\d_1,\ldots,\pm\d_\ell\}$ with $\dot\a_2\neq \pm\dot\a_3$ and $\dot\b_1,\dot\b_2\in\{\pm\ep_1,\ldots,\pm\ep_k\},$ denoting the set of nonsingular roots of $\dot R$ by $\dot R_{ns}^\times,$  we have
\[\dot\a_1+\dot\b_1=\underbrace{\underbrace{(\dot\a_2+\dot\b_2)+(\overbrace{\dot\a_3-\dot\a_2}^{\dot R_{sh}})}_{\in \dot R_{ns}^\times}+(\dot\a_1-\dot\a_3)}_{\in \dot R_{ns}^\times}+(\dot\b_1-\dot\b_2).\]
Now as each nonzero nonsingular root of $\dot R$ is of the form $\dot\a+\dot\b$ for some
$\dot\a\in\{\pm\d_1,\ldots,\pm\d_\ell\}$ and  $\dot\b\in\{\pm\ep_1,\ldots,\pm\ep_k\},$  this implies that
 for each $\dot\ep,\dot\eta\in \dot R_{ns}^\times,$
 one of the following happens:
\begin{itemize}
\item there is $\dot\b_1\in \dot R_{sh}$ such that $\dot \eta=\dot\ep+\dot \b_1,$
    \item there are   $\dot\b_1\in \dot R_{sh}$ and $\dot\b_2\in \dot R_{re}^\times$ such that
     $\dot\ep+\dot \b_1\in \dot R_{ns}^\times$ and
     $\dot \eta=\dot\ep+\dot \b_1+\dot\b_2,$
\item there are   $\dot\b_1\in \dot R_{sh}$ and $\dot\b_2,\dot\b_3\in \dot R_{re}^\times$ such that
     $\dot\ep+\dot \b_1,\dot\ep+\dot\b_1+\dot\b_2\in \dot R_{ns}^\times$ and
     $\dot \eta=\dot\ep+\dot \b_1+\dot\b_2+\dot\b_3.$
\end{itemize}
}
\end{rem}

{
\begin{deft}\label{def1}
{\rm Suppose that $S\sub R.$ We say a decomposition $S=S^+\cup S^\circ\cup S^-$ is a {\it triangular} decomposition for $S$ if there is a linear functional $\boldsymbol\zeta:\hbox{span}_\bbbr S\longrightarrow \bbbr$ such that {\small $$S^{+}=\{\a\in S\mid \boldsymbol\zeta(\a)>0\},\;S^{-}=\{\a\in S\mid \boldsymbol\zeta(\a)<0\}\andd S^\circ=\{\a\in S\mid \boldsymbol\zeta(\a)=0\}.$$}
The decomposition is called {\it trivial} if $S=S^\circ.$}
\end{deft}}

{The following proposition is crucial for the study of finite weight modules; different versions of this proposition are found in  the literature; see e.g. \cite[Pro.~3.3]{F}, \cite[\S~2]{DMP}, \cite[\S~1.4]{DG} and \cite[Pro.~2.8]{L}.}

{\begin{Pro}\label{ind}
 Suppose that $R=R^+\cup R^\circ\cup R^-$ is a nontrivial  triangular decomposition for $R$ and $R^\circ=R^{\circ,+}\cup R^{\circ,\circ}\cup R^{\circ,-}$ is a triangular decomposition for $R^\circ.$ We recall the subalgebras
$$\fl^\circ=\op_{\a\in R^{\circ,\circ}}\fl^\a,\;\; \fl^\pm=\op_{\a\in R^\pm\cup R^{\circ,\pm}}\fl^\a\andd \frak{p}=\fl^\circ\op\fl^+.$$
\begin{itemize}
\item[(i)] If $N$ is a nonzero  weight module {over} $\fl^\circ$ such that its support  $\supp(N)=\{\lam\in \fh^*\mid N^\lam\neq\{0\}\}$  lies in a single coset of $\sspan_\bbbz R^{\circ,\circ},$ then
 $$\tilde N=U(\fl)\ot_{U(\frak{p})}N$$ has a unique maximal submodule $Z$ intersecting $N$ trivially. Moreover, the induced module  $${\rm Ind}_\fl(N)=\tilde{N}/Z$$ is an irreducible   $\fl$-module if and only if $N$ is an irreducible $\fl^\circ$-module.
\item[(ii)] If $V$ is an irreducible finite weight $\fl$-module with $$V^{\fl^+}:=\{v\in V\mid \fl^+ v=\{0\}\}\neq \{0\},$$ then $V^{\fl^+}$ is an irreducible finite weight $\fl^\circ$-module and  $V\simeq {\rm Ind}_\fl(V^{\fl^+}).$
 \end{itemize}
\end{Pro}
\pf (i) As $U(\fl)$ is a free $U(\frak{p})$-module, PBW Theorem says that $\tilde N=N\op T$ in which  $T$ is an $\fh$-module. Since the support of the  $\fl^\circ$-module  $N$ is contained  in a single coset of $\sspan_\bbbz R^{\circ,\circ},$ $\supp(T)$ is disjoint from $\supp(N)$
and so $\tilde N$ contains a unique maximal submodule $Z$ intersecting $N$ trivially.}

{Next suppose that $N$ is an irreducible $\fl^\circ$-module, then each submodule of the $\fl$-module  $\tilde N$ is proper if and only if it intersects $N$ trivially and so $Z$ is the unique maximal proper submodule of $\tilde N;$ in particular, ${\rm Ind}_{\fl}(N)$ is irreducible.}

{Conversely, assume ${\rm Ind}_{\fl}(N)$ is irreducible. We know that $\fl$-module ${\rm Ind}_{\fl}(N)$ can be identified with  $N\op (T/Z)$ as an {$\fh$-module}.
If  a  nonzero weight vector $v\in T/Z$ belongs to $${\rm Ind}_{\fl}(N)^{\fl^+}=\{w\in {\rm Ind}_{\fl}(N)\mid \fl^+w=\{0\}\},$$ then  as the support of the  $\fl^\circ$-module  $N$ is contained  in a single coset of $\sspan_\bbbz R^{\circ,\circ},$ the support of the submodule generated by $v$ is disjoint from $\supp(N).$ This is a contradiction as   ${\rm Ind}_{\fl}(N)$ is irreducible. So  ${\rm Ind}_{\fl}(N)^{\fl^+}=N.$}

{Now if  $K$  is a nonzero submodule of $N,$ as above, we have ${\rm Ind}_{\fl}(K)^{\fl^+}=K.$  The assignment $\varphi:x\ot a\mapsto xa$ ($x\in U(\fl),\; a\in K$) defines an epimorphism from $U(\fl)\ot_{U(\frak{p})}K$ onto ${\rm Ind}_{\fl}(N)$ whose kernel is the unique maximal submodule intersecting $K$ trivially;
 in particular, $\varphi$ induces an isomorphism $\tilde \varphi: {\rm Ind}_{\fl}(K)^{\fl^+}\longrightarrow  {\rm Ind}_{\fl}(N)^{\fl^+}.$  Therefore, $$K=\tilde\varphi(K)=\tilde\varphi({\rm Ind}_{\fl}(K)^{\fl^+})= {\rm Ind}_{\fl}(N)^{\fl^+}=N.$$
This completes the proof.}

{(ii) Pick $0\neq v\in V^{\fl^+}.$ Then
\begin{align*}
\psi:&U(\fl)\ot_{U(\frak{p})}U(\fl^\circ) v\longrightarrow V\\
& a\ot u\mapsto au\quad\quad(a\in U(\fl),\;u\in U(\fl^\circ)v)
\end{align*}
is an epimorphism of $\fl$-modules whose kernel is the unique maximal submodule intersecting $U(\fl^\circ) v$ trivially; in particular,
$V\simeq {\rm Ind}_\fl({U(\fl^\circ) v}).$ Since $V$ is irreducible,  part (i) and its proof  implies that  $U(\fl^\circ) v$ is irreducible and ${\rm Ind}_\fl(U(\fl^\circ)v)^{\fl^+}=U(\fl^\circ)v.$ The epimorphism  $\psi$ induces an isomorphism $\tilde \psi$ from ${\rm Ind}_\fl(U(\fl^\circ)v)$ onto $V$ and we have
$$U(\fl^\circ)v=\tilde\psi(U(\fl^\circ)v)=\tilde\psi({\rm Ind}_{\fl}(U(\fl^\circ)v)^{\fl^+})= V^{\fl^+}.$$ Therefore,  $V^{\fl^+}=U(\fl^\circ)v$ is irreducible and $$V\simeq {\rm Ind}_\fl(U(\fl^\circ)v)={\rm Ind}_\fl(V^{\fl^+}).$$ This completes the proof.
\qed
}

\begin{lem}\label{trans}
Suppose that $M$ is an $\fl$-module having a weight space decomposition with respect to $\fh$ with corresponding representation $\pi.$
 Assume $ 0\neq \a\in R_{re}\cap R_0$ and  choose  $x\in\fl^\a$ and $y\in\fl^{-\a}$ such that $(x,y,h:=[x,y])$ is an $\frak{sl}_2$-triple; see (\ref{sl2}). Assume  $x$ and $y$ act locally nilpotently on $M.$ For    $\theta_\a:=\hbox{\rm exp} \pi(x)\hbox{\rm exp} \pi(-y)\hbox{\rm exp} \pi(x),$  we have $$\theta_\a(M^\lam)=M^{r_\a(\lam)}\quad\quad (\lam\in\supp(M))$$ in which  {$r_\a:\fh^*\longrightarrow \fh^*$ is defined by $r_\a(\lam):=\lam-\frac{2(\lam,\a)}{(\a,\a)}\a=\lam-\lam(h)\a$ for all $\lam\in \fh^*.$} In particular, $\lam\in \hbox{\rm supp}(M)$ if and only if  $r_\a(\lam)\in \hbox{\rm supp}(M).$
\end{lem}
\pf Since $\pi$ is a representation and $(x,y,h)$ is an $\frak{sl}_2$-triple, we have $\pi(x)=0$ if and only if $\pi(h)=0$ if and only if $\pi(y)=0.$ Also if $\pi(h)=0,$ then $\theta_\a$  as well as $r_\a\mid_{\supp(M)}$ are identity maps and so we are done. So we assume $\pi(h)\neq 0.$

Since  $(\pi(x),\pi(y),\pi(h))$ is an $\frak{sl}_2$-triple,  we have
\begin{equation}\label{exp}
\hbox{\rm exp}(\ad \pi(x))\hbox{\rm exp}(\ad (-\pi(y)))\hbox{\rm exp}(\ad \pi(x))(\pi(h))=-\pi(h).
\end{equation}

{On the other hand as} $\pi(x)$ and $\pi(y)$ are  locally nilpotent, the $\fg$-module generated by each weight vector is finite dimensional. So  the $\fg$-module $M$ is completely reducible with finite dimensional constituents and in particular, $\pi(x)$ and $\pi(y)$ are  nilpotent  on each irreducible component. We know that  if $W$ is one of these irreducible components and $T:W\longrightarrow W$ is a linear transformation,  we have
{\small $$\begin{array}{l}
 \hbox{exp}(\pi(x))~T~\hbox{exp}(-\pi(x))|_{_W}=\hbox{exp}(\ad \pi(x))(T)\andd \\\\
\hbox{exp}(\pi(-y))~T~\hbox{exp}(-\pi(-y))|_{_W}=\hbox{exp}(\ad \pi(-y))(T)
\end{array}
$$}
 and so using (\ref{exp}), we have $\theta_\a\pi(h)\theta_\a^{-1}|_{_W}=-\pi(h)|_{_W}.$ This implies that
\begin{equation}\label{theta}
\theta_\a\pi(h)\theta_\a^{-1}=-\pi(h).
\end{equation}
Now if $\lam\in\supp(M)$ and  $v\in M^\lam,$ we have  $\theta_\a(v)=\sum_{k\in\bbbz}v_{\lam+k\a}$ for some $v_{\lam+k\a}\in M^{\lam+k\a}$ $(k\in\bbbz).$ So we have
\begin{align*}
&-\lam(h)\sum_{k\in\bbbz}v_{\lam+k\a}=-\lam(h)\theta_\a(v)=-\theta_\a(\lam(h)v)\\
&=-\theta_\a(\pi(h)(v))\stackrel{(\ref{theta})}{=}\pi(h)(\theta_\a(v))=\sum_{k\in\bbbz}\pi(h)v_{\lam+k\a}=\sum_{k\in\bbbz}(\lam(h)+2k)v_{\lam+k\a}.
\end{align*}
This implies that if $v_{\lam+k\a}\neq 0$ for some $k\in\bbbz,$ then $\lam(h)+2k=-\lam(h)$ which implies that $k=-\lam(h),$ i.e., $v_{\lam+k\a}\in M^{\lam-\lam(h)\a}.$ So $\theta_\a(M^\lam)\sub M^{\lam-\lam(h)\a}=M^{r_\a(\lam)};$ similarly,  $\theta^{-1}_\a(M^{r_\a(\lam)})\sub M^{\lam}$ which  completes the proof.\qed

 \medskip
  \begin{lem}\label{trivial}
  Suppose that  $\fg$ is either $\fl$ or $\fl_0$ and $\mathcal{R}$ is the root system of $\fg$ with respect to the Cartan subalgebra  $\fh$ of $\fl,$ that is
  $$\mathcal{R}=\left\{
  \begin{array}{ll}
  R& \fg=\fl,\\
  R_0&\fg=\fl_0.
  \end{array}
  \right.  $$
  For a $\fg$-module $M$ having a weight space decomposition with respect to $\fh,$ set
 \begin{equation}\label{new2}
\hbox{\small
$\begin{array}{l}
\frak{B}_M:=\{\a\in \hbox{span}_\bbbz \mathcal{R}\mid  \{k\in\bbbz^{>0}\mid \lam+k\a\in \supp(M)\} \hbox{ is finite for all $\lam\in\supp(M)$}\} \\
\frak{C}_M:=\{\a\in \hbox{span}_\bbbz \mathcal{R}\mid   \a+\supp(M)\sub\supp(M)\}.
\end{array}$}
\end{equation}
 We also set

\begin{align}
\overline{\frak{B}}_M:=&\{\a\in {\hbox{span}_\bbbz \mathcal{R}}\mid t\a\in\frak{B}_M\hbox{ for some positive integer $t$}\},\label{saturated} \\
\overline{\frak{C}}_M:=&\{\a\in {\hbox{span}_\bbbz \mathcal{R}}\mid t\a\in\frak{C}_M\hbox{ for some positive integer $t$}\}.\nonumber
\end{align}

  We have the following:
  \begin{itemize}
  \item[(i)] Suppose $\a\in {\sspan_\bbbz \mathcal{R}}.$ Then  $\a\in \frak{B}_M$ if and only if for all  positive integers $t,$   $t\a\in \frak{B}_M$ if and only if there exists a   positive integer $t$ such that    $t\a\in \frak{B}_M;$ in particular, $\frak{B}_M=\overline{\frak{B}}_M.$
  \item[(ii)] $\a_1,\ldots,\a_n\in  \frak{C}_M$ (resp. $ \overline{\frak{C}}_M$) implies that $\a_1+\cdots+\a_n\in \frak{C}_M$ (resp. $ \overline{\frak{C}}_M$).
      \end{itemize}
  \end{lem}
  \pf (i) Suppose $\a\in \frak{B}_M$ and  $t$ is a positive integer.  As for each $\lam\in\supp(M),$ $$t\{k\in\bbbz^{>0}\mid \lam+kt\a\in \supp(M)\}\sub  \{k\in\bbbz^{>0}\mid \lam+k\a\in \supp(M)\},$$ we get that $\a\in \overline{\frak{B}}_M.$
  Next to the contrary, assume there exists a   positive integer $t\geq 2$ such that   $t\a\in \frak{B}_M$  but $\a\not\in \frak{B}_M.$ So  there is $\lam\in\supp(M)$  such that $$\aa:=\{k\in\bbbz^{>0}\mid \lam+k\a\in\supp(M)\}$$ is unbounded. Therefore, there are elements  $k_1<k_2<\cdots$ of $\aa$ and $0\leq d\leq t-1$ such that for each $i,$ $k_i\equiv d$ (mod $t$). So  $k_i=tp_i+d$ {($i\geq 1$)} for some positive integer $p_i. $ {Therefore,} we have
  \begin{align*}
  &\mu:=\lam+tp_1\a+d\a=\lam+k_1\a\in\supp(M)\andd\\
  & \mu+(p_i-p_1)t\a=\lam+tp_i\a+d\a=\lam+k_i\a\in \supp(M)
  \end{align*} for all $i\geq 2.$ This contradicts the fact that $t\a\in \frak{B}_M.$

(ii) {It is easily seen  that if $\a_1,\ldots,\a_n\in \sspan_\bbbz \mathcal{R}$ and  $t_1,\ldots,t_n\in\bbbz^{>0}$ with $t_i\a_i\in \frak{C}_M$ $(1\leq i\leq n),$ then $t_1\cdots t_n(\a_1+\cdots+\a_n)\in \frak{C}_M$.}
  \qed

  \begin{Pro}
\label{general}
{Suppose that $\fg$ is either $\fl$ or $\fl_0$ and $M$ is a $\fg$-module having a weight space decomposition with respect to $\fh$.} Denote the root system of $\fg$ with $\mathcal{R}$ and suppose that $\mathcal{S}$ is a nonempty subset of $\mathcal R$ such that
$$\hbox{$\mathcal{S}$ does not contain imaginary roots},\;\;\mathcal{S}\sub \frak{B}_M\andd -\mathcal{S}\sub \frak{C}_M.$$
Then we have the following:
\begin{itemize}
\item[(i)] If $\aA$ is  a nonempty subset of $\supp(M)$ with $(\aA+\mathcal{S})\cap \supp(M)\sub \aA,$ then for each $\b\in \mathcal{S},$ $$\aA_\b:=\{\lam\in\aA\mid \lam+\b\not\in \supp(M)\}$$ is also  nonempty with $(\aA_\b+\mathcal{S})\cap \supp(M)\sub \aA_\b.$
\item[(ii)]  If $\mathcal{S}$ is finite and  $\aA$ is as in part {\rm (i)}, then there is $\lam\in \aA$ such that $$(\lam+\sspan_{\mathbb{Z}^{\geq 0}}\mathcal{S})\cap \supp(M)=\{\lam\}.$$
\end{itemize}
\end{Pro}
\pf (i) Suppose that $\lam\in \aA$ and $\b\in \mathcal{S}.$ Since $\b\in \frak{B}_M,$ there is a nonnegative integer $k$ such that $\mu:=\lam+k\b\in \supp(M)$ and $\mu+\b\not\in \supp(M).$ We claim that $\mu\in \aA_\b.$ We just need to show $\mu\in\aA.$ Since $-\b\in \frak{C}_M,$  $\lam+(k-t)\b\in \supp(M)$ for all $0\leq t\leq k.$  Since $(\aA+\mathcal{S})\cap \supp(M)\sub \aA,$ it follows that $\lam+(k-t)\b\in \aA$ for all $0\leq t\leq k;$  in particular, $\mu\in \aA.$

To complete the proof, we need to show $(\aA_\b+\mathcal{S})\cap \supp(M)\sub \aA_\b.$ Suppose $\nu\in\aA_\b$ and $\gamma\in \mathcal{S}$ are such that $\nu+\gamma\in\supp(M).$ If {to the contrary,} $\nu+\gamma+\b\in \supp(M),$ since $-\gamma\in \frak{C}_M,$ we get $\nu+\b\in \supp(M)$ which contradicts the fact that $\nu\in \aA_\b.$ So $\nu+\gamma+\b\not\in \supp(M);$ in other words, $\nu+\gamma\in \aA_\b.$

(ii) Suppose $\mathcal{S}=\{\b_1,\ldots,\b_N\}.$
Set $$\aA_0:=\aA,\;\; \aA_{t+1}:=(\aA_{t})_{\b_{t+1}}=\{\lam\in\aA_t\mid \lam+\b_{t+1}\not\in \supp(M) \}\;\;\;(0\leq t\leq N-1).$$
We have $\aA_N\sub \aA_{N-1}\sub\cdots\sub \aA_1$ and by part (i), for each $1\leq t\leq N,$ $\aA_t\neq \emptyset;$ in particular, $\aA_N\neq \emptyset.$ For $\lam\in \aA_N,$ since $\lam\in \aA_t$ $(1\leq t\leq N),$ we get $\lam+\b_t\not \in \supp(M)$ which in turn implies that $(\lam+\sspan_{\mathbb{Z}^{\geq 0}}\mathcal{S})\cap \supp(M)=\{\lam\}$ as $-\mathcal{S}\sub \frak{C}_M.$
\qed

\begin{Pro}\label{pdelta}
  {Suppose that $\fg$ is either $\fl$ or $\fl_0$ and $M$ is a $\fg$-module having a weight space decomposition with respect to $\fh.$  Denote the root system of $\fg$ with respect to  $\fh$  with $\mathcal{R}.$ Assume     $\mathcal{R}=\mathcal{R}^+\cup \mathcal{R}^\circ\cup \mathcal{R}^-$ is  a triangular decomposition for $\mathcal{R}$ with corresponding functional $\boldsymbol\zeta.$ Set}
     $$\mathcal{R}_{re}^{\pm}:=\mathcal{R}^\pm\cap \mathcal{R}_{re}\andd \mathcal{R}_{im}^\pm:=\mathcal{R}_{im}\cap \mathcal{R}^\pm.$$
     Assume $\mathcal{R}_{re}^+ \sub \frak{B}_M,$ $\mathcal{R}^-_{re}\sub \frak{C}_M;$ see (\ref{new2})  and $\boldsymbol\zeta(\d)>0.$
   If  $p\in\bbbz^{>0}$ and $\lam\in \supp(M)$ are  such that
    $(\lam+\bbbz^{>0}p\d)\cap\supp(M)=\emptyset,$ then there is $\mu\in\supp(M)$ such that $(\mu+(\mathcal{R}_{re}^+\cup \mathcal{R}_{im}^+))\cap\supp(M)=\emptyset.$
\end{Pro}
\pf
Set
\[\dot{\mathcal{R}}:=\left\{
\begin{array}{ll}
\dot R& \fg=\fl\\
\{\dot \a\in \dot R\mid R_0\cap (\dot\a+\bbbz\d)\neq \emptyset\}& \fg=\fl_0.
\end{array}
\right.
\]
  Using  (\ref{imp}), one knows that for each $0\neq \dot\a\in \dot {\mathcal{R}},$
there is $r_{\dot\a}\in \bbbz^{>0}$ and $k_{\dot\a}\in \bbbz^{\geq 0}$ such that
\begin{equation}\label{ralpha}
\{n\in\bbbz\mid \dot\a+n\d\in \mathcal{R}\}=r_{\dot\a}\bbbz+k_{\dot\a}
\end{equation}
and that
\begin{equation}\label{new}
\parbox{3.6in}{there is $0\neq \dot\a^*\in \dot{\mathcal{R}}_{re}$ such that $k_{\dot\a^*}=0$ and  $r_{\dot\a^*}\bbbz\d=\mathcal{R}_{im}.$}
\end{equation}
Fix $\lam$ and $p$ as in the statement. Consider (\ref{ralpha}) and for $\dot\a\in \dot{\mathcal{R}}_{re}^\times,$ suppose that
\begin{equation}\label{tdotalpha}
\parbox{4 in}{
$t_{\dot \a}\in r_{\dot\a}\bbbz$ is the smallest integer  such that ${\boldsymbol{\zeta}}(\dot\a+(t_{\dot\a}+k_{\dot\a})\d)>0.$}
\end{equation}
  Set
\begin{align}
\mathcal{P}:=&\{\dot\a+(t_{\dot\a}+k_{\dot\a}+s)\d\mid \dot\a\in \dot{\mathcal{R}}_{re}^\times,\;0\leq s\leq r_{\dot\a}p\}\cap \mathcal{R}\sub \mathcal{R}_{re}^+,\nonumber\\
\label{s}
\mathcal{S}:=&\{\dot\a+(t_{\dot\a}+k_{\dot\a})\d\mid \dot\a\in \dot{\mathcal{R}}_{re}^\times\}\sub \mathcal{P}
\end{align}
and

$$\aa:=\{\mu\in {\rm supp}(M)\mid \{\a\in \mathcal{R}_{re}^+\mid \mu+\a\in{\rm supp}(M)\}\sub\mathcal{P}\}.$$
We {have in particular that}
\begin{equation}\label{finite}
\parbox{4in}{
if $\mu\in \aa,$ then $\{\a\in \mathcal{R}_{re}^+\mid \mu+\a\in{\rm supp}(M)\}$ is a finite set.
}
\end{equation}

\smallskip

\noindent {\bf Claim 1.} $\aa$ is a nonempty set:
 We claim that $\lam$ as in the statement belongs to $\aa.$
Suppose $\a\in \mathcal{R}_{re}^+$ is such that $\lam+\a\in\supp(M).$ We shall show $\a\in\mathcal{P}.$ Since $\a\in \mathcal{R}_{re}^+,$ by (\ref{ralpha}) and (\ref{tdotalpha}),
$$\a=\dot\a+m\d+k_{\dot\a}\d\quad\hbox{for some $\dot\a\in \dot{\mathcal{R}}_{re}^\times$ and $m\in r_{\dot\a}\bbbz$ with $m\geq t_{\dot\a}.$}$$
We have $m-t_{\dot\a}=kr_{\dot\a}p+s$ for some nonnegative integer $k$ and $s\in \{0,\ldots, r_{\dot\a}p\}.$ We notice that as $r_{\dot\a}| m$ and $r_{\dot\a}|t_{\dot\a} ,$  we have $r_{\dot\a}| s;$ in particular, $\dot\a+(t_{\dot\a}+k_{\dot\a}+s)\d\in \mathcal{R}_{re}^+.$
We also have
\begin{align*}
\lam+\dot\a+(t_{\dot\a}+k_{\dot\a}+s)\d+kr_{\dot\a}p\d=&\lam+\dot\a+(t_{\dot\a}+k_{\dot\a})\d+(m-t_{\dot\a})\d\\
=&\lam+\dot\a+k_{\dot\a}\d+(t_{\dot\a}\d+(m-t_{\dot\a})\d)\\
=&
\lam+(\dot\a+(m+k_{\dot\a})\d)=\lam+\a\in\supp(M).
\end{align*}
Since   $-(\dot\a+(t_{\dot\a}+k_{\dot\a}+s)\d)\in \mathcal{R}_{re}^-\sub \frak{C}_M,$ we conclude $\lam+kr_{\dot\a}p\d\in\supp(M)$ which implies that $k=0$ by our assumption on {$p$ and $\lam.$} So $\a=\dot\a+(t_{\dot\a}+k_{\dot\a}+s)\d\in \mathcal{P}.$

\smallskip

\noindent{\bf Claim 2.} For each $\mu\in\aa,$ $\{m\d\in\mathcal{R}_{im}^+\mid \mu+m\d\in {\rm supp}(M)\}$ is a finite set:
Suppose  $\mu\in\aa$ and to the contrary assume
\begin{center}
there are infinitely many $m\d\in\mathcal{R}_{im}^+$ such that $\mu+m\d\in \hbox{supp}(M).$
\end{center}
We know from (\ref{new}) and Table~\ref{table2} that there is $\dot\a^*\in \dot{\mathcal{R}}_{re}^\times $ such that
$$\{n{\d}\in\bbbz\mid -\dot\a^*+n\d\in \mathcal{R}\}=\{n{\d}\in\bbbz\mid \dot\a^*+n\d\in \mathcal{R}\}=r_{\dot\a^*}\bbbz{\d}=\mathcal{R}_{im}$$ and
\begin{equation}
\label{new3}
r_{\dot\a^*}\mid r_{\dot\a}\quad\quad(\dot\a\in \dot R^\times).
\end{equation}
So there are infinitely many $m\in r_{\dot\a^*}\bbbz$ such that $m\geq t_{\dot\a^*}$ (see (\ref{tdotalpha})) and $\mu+m\d\in \supp(M).$
Since $-(\dot\a^*+t_{\dot\a^*}\d)\in \mathcal{R}^-\sub \frak{C}_M,$ we get that $\mu+(-\dot\a^*+(m-t_{\dot\a^*})\d)\in\supp(M)$ for infinitely many $m\in r_{\dot\a^*}\bbbz$ with $m>t_{\dot\a^*}.$ But this contradicts (\ref{finite}) as $\mu\in\aa.$
\smallskip

\noindent{\bf Claim 3.} There is $\mu\in \supp(M)$ such that $\mu+m\d\not\in \supp(M)$ for all $m\d\in \mathcal{R}^+_{im}:$ Pick   $\eta\in\aa.$ Using Claim 2, we assume  $N$ is the greatest nonnegative  integer of $r_{\dot\a^*}\bbbz$ with $\eta+N\d\in \supp(M).$ So for $\mu:=\eta+N\d$ and $m\d\in r_{\dot\a^*}\bbbz^{\gneq 0}\d=\mathcal{R}_{im}^+,$ $\mu+m\d\not\in \supp(M).$

\smallskip

\noindent{\bf Claim 4.}  Set $ X:=\{\mu\in\supp(M)\mid \forall m\d\in\mathcal{R}^+_{im},\;\mu+m\d\not\in \supp(M)\}.$ Recall (\ref{s}), then there is $\mu\in X$ such that $(\mu+\sspan_{\bbbz^{\geq0}}\mathcal{S})\cap \supp(M)=\{\mu\}:$
Using Proposition \ref{general}(ii) and Claim 3, we need to show $(X+\mathcal{S})\cap\supp(M)\sub X.$ To the contrary assume   $\mu\in X$ and $\b\in \mathcal{S}$ are  such that $\mu+\b\in \supp(M)$ and $\mu+\b\not\in X.$ So there is  $m\d\in \mathcal{R}_{im}^+=r_{\dot\a^*}\bbbz^{\gneq0}\d$ such that  $\mu+\b+m\d\in\supp(M),$ then as $-\b\in \frak{C}_M,$ $\mu+m\d\in \supp(M)$ which is a contradiction as $\mu\in X$.

\smallskip

\noindent{\bf Claim 5.} There is $\mu\in \hbox{supp}(M)$ such that $(\mu+(\mathcal{R}_{re}^+\cup \mathcal{R}_{im}^+))\cap \hbox{supp}(M)=\emptyset:$ Using Claim 4, we  choose $\mu\in \supp(M)$ such that $$(\mu+(\mathcal{R}^+_{im}\cup \sspan_{\bbbz^{\geq0}}\mathcal{S}))\cap \supp(M)=\{\mu\}$$ If  $\a\in \mathcal{R}_{re}^+\cup \mathcal{R}_{im}^+$ and $\mu+\a\in \supp(M),$ then  $\a\in \mathcal{R}_{re}^+.$ So $\a=\dot\a+m\d+k_{\dot\a}\d$ for some  $\dot\a\in \dot{\mathcal{R}}_{re}^\times$  and some integer $m\in r_{\dot\a}\bbbz$ with $m\geq t_{\dot\a};$ see (\ref{tdotalpha}).  If $m\gneq t_{\dot\a},$  we get $\mu+(m-t_{\dot\a})\d=\mu+\a-(\dot\a+k_{\dot\a}\d+t_{\dot\a}\d)\in \supp(M)$ as  $-(\dot\a+t_{\dot\a}\d+k_{\dot\a}\d)\in \mathcal{R}_{re}^{-}\sub \frak{C}_M, $ and  $\mu+\a\in \supp(M).$ But this  contradicts the choice of $\mu$ as by (\ref{new3}), $$(m-t_{\dot\a})\d\in r_{\dot\a}\bbbz^{\gneq0}\d\sub r_{\dot\a^*}\bbbz^{\gneq0}\d=\mathcal{R}^+_{im}.$$ So $m=t_{\dot\a};$ i.e., $\a\in \mathcal{S}.$ It means that $$\mu\neq\mu+(\dot\a+t_{\dot\a}\d+k_{\dot\a}\d)=\mu+\a\in \supp(M)\cap (\mu+\mathcal{S})$$ which is  again a contradiction. So there is no $\a\in \mathcal{R}_{re}^+\cup \mathcal{R}_{im}^+$ with  $\mu+\a\in \supp(M).$ This completes the proof.
\qed

\begin{Pro}\label{last-one}
Recall $\dot R$ from Table~\ref{dot} and assume $M$ is a module over the affine Lie superalgebra  $\fl.$ Suppose $\boldsymbol\zeta$ is a linear functional on $\sspan_\bbbr R$ with corresponding triangular decomposition $R=R^+\cup R^\circ\cup R^-.$ Set \begin{align*}
\aa:=&\{v\in M\setminus\{0\}\mid \fl^{\a}v= \{0\};\;\; \;\;\forall\a\in  R^+\cap( R_{re}\cup R_{im})\}\\
=&\{v\in M\setminus\{0\}\mid \fl^{n\d}v= \fl^{\a}v= \{0\};\;\; \;\;\forall\a\in  R_{re}\cap R^+,\; n\in\bbbz^{>0}\},
\end{align*}
and  assume
$$B:=\{v\in\aa\mid \forall\dot \a\in \dot R_{ns}^\times~ \exists~ N\in\bbbz^{\geq0}\hbox{ s.t. } \fl^{\dot\a+n\d} v=\{0\}\;\;\;(\forall  n\geq N)\}$$ is nonempty. If $\boldsymbol\zeta(\d)>0,$ then
$$M^{\fl^+}=\{v\in M\mid \fl^\a v=\{0\}\;\; (\forall \a\in R\hbox{ with } \boldsymbol\zeta(\a)>0) \}\neq \{0\}.$$
\end{Pro}
\pf We know from (\ref{imp}) and Table~\ref{table2} that  for each $\dot\a\in \dot R^\times,$ there is $r_{\dot\a}\in\bbbz^{>0}$  and $0\leq k_{\dot\a}<r_{\dot\a}$ such that
\begin{equation}\label{salphadot1}
\begin{array}{l}
S_{\dot\a}=\{m\d\mid m\in \bbbz,\;\;\dot\a+m\d\in R\}=(r_{\dot\a}\bbbz+k_{\dot\a})\d;\quad(\dot\a\in \dot R^\times),\\
k_{\dot\a}=k_{\dot\b}=0,\;r_{\dot\b}=r_{\dot\a};\quad(\dot\a,\dot\b\in \dot R_{ns}^\times).
\end{array}
\end{equation}
In particular,
\begin{equation}\label{ss}
\parbox{2in}{$S_{\dot\a}$ is a group for all $\dot\a\in \dot R^\times_{ns}.$}
\end{equation}
Since ${\boldsymbol{\zeta}}(\d)>0,$ for each $0\neq \dot\a\in \dot R,$  we assume
\begin{equation}\label{m1}
\hbox{$m_{\dot\a}$ is the smallest integer such that for  $\b_{\dot\a}:=\dot\a+(r_{\dot\a}m_{\dot\a}+k_{\dot\a})\d\in R,$ $\boldsymbol\zeta(\b_{\dot\a})>0.$}
\end{equation}
Set
$$ \Phi:=\{\b_{\dot\a}\mid \dot\a\in \dot R^\times\}.$$

\noindent{\bf Claim 1.} {\small $B=B':=\{v\in\aa\mid \exists N\in\bbbz^{\geq0}~s.t.~ \fl^{\a+n\d}v=\{0\}\quad(\a\in \Phi\cap R_{ns},\;n\geq N)\}:$}
Suppose that $v\in B.$ So for each $\dot\a\in \dot R_{ns}^\times,$ there is $N_{\dot\a}\in\bbbz^{\geq 0}$ with $\fl^{\dot\a+n\d}v=\{0\}$ for all $n\geq N_{\dot\a}.$ Set $N:=max\{N_{\dot\a}-(r_{\dot\a} m_{\dot\a}+k_{\dot\a})\mid \dot\a\in \dot R_{ns}^\times\}.$ Then
$\fl^{\b_{\dot\a}+n\d}v=\{0\}$ for all $n\geq N$ and $\dot\a\in \dot R_{ns}^\times,$ i.e., $B\sub B'.$ Conversely,
suppose $v\in B'$ and pick $N\in\bbbz^{\geq0}$ with $\fl^{\b_{\dot\a}+n\d}v=\{0\}$ for $\dot\a\in \dot R^\times_{ns}$ and $n\geq N .$ So for each $\dot\a\in \dot R_{ns}^\times$ and $n\geq N+(r_{\dot\a}m_{\dot\a}+k_{\dot\a}),$
 we have $\fl^{\dot\a+n\d}v=\{0\},$ that is $v\in B.$

\medskip
Using  Claim~1, for $v\in B,$  we set
$$n_v:=min\{ N\in\bbbz^{\geq 0}\mid \fl^{\a+n\d} v=\{0\}\;\;\;(\a\in \Phi\cap  R_{ns},\;\;  n\geq N)\}$$and
 $$ C_v:=\{\a+t\d\mid \a\in \Phi\cap  R_{ns},\; 0\leq t< n_v\}\cap R\sub R_{ns}.$$

\noindent{\bf Claim 2.}
Assume $v\in B,$    $N\in \bbbz^{\geq 0}$ and $\a\in  C_v$ satisfy
\begin{itemize}
\item[(1)] $\fl^{\a+N\d}v\neq \{0\},$
\item[(2)] if $\a'\in  C_v$ and $\fl^{\a'+N\d}v\neq 0,$ then $\boldsymbol\zeta(\a')\leq \boldsymbol\zeta (\a),$
\item[(3)] for all positive integers $m$ and $\a'\in  C_v,$ $\fl^{\a'+N\d+m\d}v=\{0\}.$
\end{itemize}
Then for  $0\neq w\in \fl^{\a+N\d}v,$ $w\in B:$
We carry out this in the following stages:

\noindent{\underline{Stage 1}.} For $m\in\bbbz^{>0},$ $\fl^{m\d}w=\{0\}:$
Use (3) and note that $v\in \aa$ to get that  $$\fl^{m\d}w\sub \fl^{m\d}\fl^{\a+N\d}v\sub \underbrace{\fl^{\a+(N+m)\d}v}_{0}+\fl^{\a+N\d}\underbrace{\fl^{m\d}v}_{0}=\{0\}.$$

\noindent{\underline{Stage 2}.}  For $\b\in R_{re}^\times$ with $\boldsymbol\zeta(\b)>0,$
$\fl^\b w=\{0\}:$ Since $v\in\aa,$ $\fl^\b v=\{0\},$ so
we have
\begin{equation}\label{new4}\fl^{\b}w\sub \fl^{\b}\fl^{\a+N\d}v\sub \fl^{\a+\b+N\d}v+\fl^{\a+N\d}\fl^{\b}v=\fl^{\a+\b+N\d}v .
\end{equation}
The following three cases can happen:
\begin{itemize}
\item    $\a+\b+N\d\not\in R:$ Then $\fl^{\b}w\sub \fl^{\a+\b+N\d}v=\{0\}.$
\item  $\a+\b+N\d\in  R_{re}:$ As $v\in \aa$ and $\boldsymbol\zeta(\a+\b+N\d)=\underbrace{\boldsymbol\zeta(\a)}_{>0}+\underbrace{\boldsymbol\zeta(\b)}_{>0}+\underbrace{\boldsymbol\zeta(N\d)}_{\geq0}>0,$ we get that  $\fl^{\a+\b+N\d}v=\{0\}$ and so $\fl^\b w=\{0\}$.
    \item     $\a+\b+N\d\in R_{ns}^\times:$
Regarding (\ref{salphadot1}), suppose $\a=\dot\a+\sg$ and $\b=\dot\b+\tau$ for some $\dot\a,\dot\b\in \dot R^\times,$ $\sg\in S_{\dot\a}$ and $\tau\in S_{\dot\b}.$ Since $\a+\b+N\d\in R^\times_{ns},$ $\dot\gamma:=\dot\a+\dot\b\in \dot R^\times_{ns}.$ So we have
\begin{align*}
\hbox{\small $\left\{
\begin{array}{l}
\dot\a+\sg=\a\in R^\times_{ns}\\
\dot\a+\sg+N\d=\a+N\d\in R^\times_{ns}\\
\dot\gamma+\sg+\tau+N\d=\a+\b+N\d\in  R^\times_{ns}
\end{array}
\right. $}&\hbox{\small $\longrightarrow\left\{
\begin{array}{l}
\sg\in S_{\dot\a}\\
\sg+N\d\in S_{\dot\a}\\
\sg+\tau+N\d\in S_{\dot\gamma}
\end{array}
\right.$}\\
&\hbox{\small$\xrightarrow{(\ref{salphadot1}),(\ref{ss})}$}\hbox{\small$\sg+\tau\in S_{\dot\gamma}.$}
\end{align*} So $\a+\b=\dot\gamma+(\sg+\tau)\in R^\times_{ns}.$ Since  $\boldsymbol\zeta(\a+\b)>0,$ by (\ref{m1}), there exists $m'\in\bbbz^{\geq0}$ such that  $$\a+\b=\gamma+m'\d\quad\hbox{where}\quad \gamma:=\dot\gamma+(r_{\dot\gamma}m_{\dot\gamma}+k_{\dot\gamma})\d\in \Phi\cap R_{ns}\sub C_v.$$ So
\begin{align*}
\a+\b+N\d
=\gamma+(m'+N)\d.
\end{align*}
If  $m'=0,$ then $\a+\b=\gamma\in C_v$ and as $\boldsymbol\zeta(\gamma)=\boldsymbol\zeta(\a+\b)>\boldsymbol\zeta(\a),$ using (2), we have  $$\fl^{\b}w\stackrel{(\ref{new4})}{\sub}\fl^{\a+\b+N\d}v=\fl^{\gamma+N\d}v=\{0\}.$$
Also if  $m'>0,$ then   (3) implies that  $$\fl^{\b}w\stackrel{(\ref{new4})}\sub\fl^{\a+\b+N\d}v=\fl^{\gamma+(m'+N)\d}v=\{0\}.$$
\end{itemize}

\noindent{\underline{Stage 3.}} $w\in B:$ Contemplating Claim~1 and using Stages~1,2, we need to show that there is a positive integer $P$ such that  for all $\eta\in\Phi \cap R_{ns}$ and   $n\geq P,$ $\fl^{\eta+n\d} w=\{0\}$.
Since $v\in B,$ we pick $P\in\bbbz^{>0}$ such that $\fl^{\eta+n\d} v=\{0\}$ for all $\eta\in\Phi \cap R_{ns}$ and   $n\geq P.$ Then for all $\eta\in \Phi\cap R_{ns}$ and $n\geq P,$ we have $$\fl^{\eta+n\d}w\sub\fl^{\eta+n\d}\fl^{\a+N\d}v\sub\fl^{\eta+\a+n\d+N\d}v+\fl^{\a+N\d}\underbrace{\fl^{\eta+n\d}v}_0=\fl^{\eta+\a+n\d+N\d}v.$$ But if $\eta+\a+n\d+N\d\in R,$ then by (\ref{re-re}), $\eta+\a+n\d+N\d\in R_{re},$ so as $v\in\aa$ and $\boldsymbol\zeta(\eta+\a+n\d+N\d)>0,$ we get
$\fl^{\eta+\a+n\d+N\d}v=\{0\}.$ Therefore, we have  $\fl^{\eta+n\d}w=\{0\}.$

\noindent{\bf Claim  3.}
For $v\in B,$ $n_v\neq 0$ if and only if $$\aa_v:=\{\a\in C_v\mid  \fl^{\a+m\d}v\neq \{0\}\hbox{ for some $m\geq 0$}\}\sub R_{ns}$$ is a nonempty set:
It follows from the following: \begin{align*}
n_{v}=0\Leftrightarrow&\fl^{\a+m\d}v=\{0\}\quad\quad (\a\in \Phi\cap R_{ns},\; m\geq0)\\
\Leftrightarrow&\fl^{\a+m\d+t\d}v=\{0\}\quad\quad (\a\in \Phi\cap R_{ns},\; m\geq0,\; 0\leq t\leq n_{v})\\
\Leftrightarrow&\fl^{\a+m\d}v=\{0\}\quad\quad (\a\in  C_{v},\; m\geq0)\\
\Leftrightarrow&\aa_{v}=\emptyset.
\end{align*}

\noindent{\bf Claim 4.} If $v\in B$ and $\aa_v\neq\emptyset,$ then there is $0\leq k<n_{v}$ such that $$B_k(v):=\{\a\in C_{v}\mid \fl^{\a+k\d}v\neq \{0\}\}$$ is nonempty: Since $\aa_v\neq \emptyset,$ there is $\a\in C_v$ and $k\in\bbbz^{\geq 0}$ such that
$\fl^{\a+k\d}v\neq \{0\}.$ Since $\a\in C_v,$ there is $\b\in \Phi\cap R_{ns}$ and $0\leq t<n_v$ such that $\a=\b+t\d.$ So
$\{0\}\neq \fl^{\b+(k+t)\d}v.$ Therefore, we get $0\leq k\leq k+t<n_v.$

\noindent{\bf Claim 5.} For $v\in B$ with $\aa_v\neq\emptyset,$  set  $$N(v):={\rm max}\{0\leq k<n_{v}\mid B_k(v)\neq \emptyset\}$$ where $B_k(v)$ is as in the previous claim and
 choose $\ep\in B_{N(v)}(v)$ with $$\boldsymbol\zeta(\ep)={\rm max}\{\boldsymbol\zeta(\a)\mid \a\in B_{N(v)}(v)\}.$$ Then for   $0\neq w\in\fl^{\ep+N(v)\d}v,$ $w\in B$ and $\ep+N(v)\d\in \aa_{v}\setminus\aa_w:$ That $w\in B$ follows from Claim~2.
We shall show  $\ep+N(v)\d\in \aa_{v}\setminus\aa_w.$

 Since $\ep\in C_{v},$ there is $\eta\in \Phi\cap R_{ns}$ and   $1\leq p<n_{v}$ with $\ep=\eta+p\d.$ But $\fl^{\eta+(p+N(v))\d}v=\fl^{\ep+N(v)\d}v\neq \{0\},$ so $p+N(v)< n_{v},$ in other words, $$\ep+N(v)\d=\eta+(p+N(v))\d\in C_{v}$$ and $\fl^{\ep+N(v)\d}v\neq \{0\}$ which means that $\ep+N(v)\d\in \aa_{v}.$ So, we just need to show $\ep+N(v)\d\not\in \aa_w.$
  Since $N(v)={\rm max}\{0\leq k<n_{v}\mid B_k(v)\neq \emptyset\},$  we have $$\fl^{\ep+N(v)\d+n\d}v=\{0\}\quad(n>0).$$  This together with the fact that two times of  a nonzero nonsingular root is not a root, gives that
\begin{align*}
\fl^{\ep+N(v)\d+n\d}w\sub \fl^{\ep+N(v)\d+n\d}\fl^{\ep+N(v)\d}v=&\fl^{\ep+N(v)\d}\fl^{\ep+N(v)\d+n\d}v=\{0\}\;\; (n>0)\andd \\
\fl^{\ep+N(v)\d}w\sub \fl^{\ep+N(v)\d}\fl^{\ep+N(v)\d}v=&[\fl^{\ep+N(v)\d},\fl^{\ep+N(v)\d}]v=\{0\}.
\end{align*} Therefore, $\ep+N(v)\d\not\in \aa_w$ as we desired.

\noindent{\bf Claim 6.} There is $v_0\in B$ such that $n_{v_0}=0,$ i.e., $v_0\in M^{\fl^+}:$ Assume $v_0\in B$ is such that\footnote{We use $|X|$ to denote the cardinal number of a set  $X.$}   $$|\aa_{v_0}|=\hbox{\rm min}\{|\aa_v|\mid v\in B\}.$$ We  claim that $n_{v_0}=0.$
To the contrary, assume $n_{v_0}\neq 0.$ By Claim~3, $\aa_{v_0}\neq \emptyset.$ Choose $\ep$ and $N(v_0)$ as in Claim~5 and  pick a nonzero element  $w\in\fl^{\ep+N(v_0)\d}v_0.$ So by Claim~5,  $w\in B.$
If $\a\in \aa_w,$ then there is $m\in\bbbz^{\geq 0}$ such that
$$\{0\}\neq \fl^{\a+m\d}w\sub \fl^{\a+m\d+\ep+N(v_0)\d}v_0+\fl^{\ep+N(v_0)\d}\fl^{\a+m\d}v_0.$$
But either  $\a+m\d+\ep+N(v_0)\d\not \in R$ or $\a+m\d+\ep+N(v_0)\d\in R_{re}\cup  R_{im}$  (see (\ref{re-re})) with $\boldsymbol\zeta(\a+m\d+\ep+N(v_0)\d)>0,$ so $\fl^{\a+m\d+\ep+N(v_0)\d}v_0=\{0\},$ i.e., $$\{0\}\neq \fl^{\a+m\d}w\sub\fl^{\ep+N(v_0)\d}\fl^{\a+m\d}v_0$$ which in turn implies that  $\fl^{\a+m\d}v_0\neq\{0\},$ that is, $\a\in \aa_{v_0}.$ This means that $$\aa_w\subseteq \aa_{v_0}.$$ But by Claim~5, $ \aa_{v_0}\setminus\aa_w\neq \emptyset$ which is a contradiction as $\aa_{v_0}$ has the minimum  cardinality among all $\aa_u$ ($u\in B$).
\qed
\section{Modules having shadow}
Keep the same {notations} as in Section \ref{generic} {and assume $M$ is a weight $\fl$-module.}
  Denote by $R^{in}$ (resp. $R^{ln}$) the set of all nonzero $\a\in R_{re}$ for which  $0\neq x\in \fl^\a$ acts injectively (resp. locally nilpotently) on $M$.

\begin{deft}\label{shadow}
 {\rm   We say  $M$ has {\it shadow} if
 \begin{itemize}
 \item[\rm\bf (s1)]
$R_{re}^\times=R^{in}\cup R^{ln},$
\item[\rm\bf (s2)] $R^{ln}=\frak{B}_M\cap R_{re}^\times$ and $R^{in}=\frak{C}_M\cap R_{re}^\times.$
 \end{itemize}}
\end{deft}
\begin{rem}
{\rm We mention that if the $\fl$-module $M$ has shadow, then  $\a\in R^{ln}$ (resp. $\a\in R^{in}$) if and only if
  $\{k\in\bbbz^{\geq 0}\mid\lam+k\a\in\supp(M)\}$ is bounded (resp. unbounded) for some $\lam\in \supp(M).$}
\end{rem}
\begin{lem}\label{nilpot}
Suppose that $\gG$ is a Lie superalgebra and $\phi:\gG\longrightarrow \End{V}$ is a representation of  $\gG$ in a superspace $V.$ For each nonnegative integer $n,$ define
\[b_{2i}^n:=b_{2i}^{n-1}+b_{2i-2}^{n-1}\;\;\;\; (n\geq2,\;1\leq i\leq n-1)\andd b_0^n=b_{2n}^n:=1.\] Then for $n\in\bbbz^{\geq 0}$ and homogeneous elements $x,y\in\gG,$ if $|y|=1,$  we have
\begin{align*}
\phi(y)^{2n}\phi(x)&= \sum_{i=0}^n b_{2i}^n\phi({\rm ad}y^{2i}(x))\phi(y)^{2n-2i}\andd\\
\phi(y)^{2n+1}\phi(x)&= \sum_{i=0}^n b_{2i}^n((-1)^{|x|}\phi({\rm ad}y^{2i}(x))\phi(y)^{2n+1-2i}+\phi({\rm ad}y^{2i+1}(x))\phi(y)^{2n-2i})\\
\end{align*}
and if  $|y|=0,$ we have
\[\phi(y)^n\phi(x)=\sum_{i=0}^n\binom{n}{i}\phi(({\rm ad}y)^i(x))\phi(y)^{n-i}\quad\quad\quad(n\in\bbbz^{\geq0}).\]
\end{lem}
\pf It is easily verified.\qed

\begin{Pro}\label{suff}
\begin{itemize}
\item[(i)]  Suppose that the $\fl$-module $M$ is irreducible, then {\bf (s1)} is satisfied.
\item[(ii)] Suppose that the  $\fl$-module $M$ satisfying {\bf (s1)} and each weight space is finite dimensional. Then $M$ has shadow.
\end{itemize}
\end{Pro}
\pf
(i) It follows from Lemma \ref{nilpot}.

(ii)
It is trivial that if $\a\in R^{in},$ then $\a\in \frak{C}_M,$ so to complete the proof, we just need to assume $\a\in R^{ln}$ and show that $\{k\in\bbbz^{\geq 0}\mid \lam+k\a\in\supp(M)\}$ is bounded for all $\lam\in\supp(M).$ Two cases can happen: $-\a\in R^{ln}$ and $-\a\in R^{in}.$
We need to study  separately each case for $\a\in R_1$ and $\a\in R_0.$

We first study the case that  $\a\in R^{ln}$ is a real odd root.
Fix  $x\in\LL^\a$ and  $y\in\LL^{-\a}$ such that $$\fg:=\hbox{span}_\bbbc\{x,y,h:=[x,y],[x,x],[y,y]\}$$ is a Lie superalgebra isomorphic to $\frak{osp}(1,2)$ with $\a(h)=2;$ see \cite[\S~3]{you3} and \cite[Exa. 2.2]{AY}.

To get the result in this case, we first assume $-\a\in R^{ln}.$
For each $\lam\in\supp(M),$ $W:=\op_{k\in\bbbz} M^{\lam+k\a}$ is a $\fg$-module. The set of eigenvalues of the action of $h$ on $W:=\op_{k\in\bbbz} M^{\lam+k\a}$ is $\Lam:=\{\lam(h)+2k\mid k\in\bbbz,\; \lam+k\a\in\supp(M)\}$ and the eigenspace corresponding to  each $\lam(h)+2k\in\Lam$ is the finite dimensional space $M^{\lam+k\a}.$

Since both $x$ and $y$ act locally  nilpotently, the $\fg$-submodule of $W$  generated by a weight vector is finite dimensional. So  it follows from  \cite[Thm. 2.6]{you3} that $W$ is  completely reducible with finite dimensional irreducible constituents.  In particular, by \cite[Lem.~2.4(iii)]{you3},  dimension of the  eigenspace corresponding to $0$ is infinite if there are infinitely  many constituents. But the  eigenspace corresponding to $0$ is $M^{\lam-(\lam(h)/2)\a}$ which is finite dimensional. Therefore, there are just finitely many constituents and so again using \cite[Lem.~2.4(iii)]{you3},
 $\{k\in\bbbz\mid\lam+k\a\in\supp(M)\}$ is bounded and so we are done in the case that $\pm\a\in R^{ln}\cap R_1.$

Next assume $\a\in R^{ln}\cap R_1$ and $-\a\in R^{in}.$ For a positive integer $m$ and a weight $\nu,$ set
$$
r_m(\nu)=\left\{
\begin{array}{ll}
\displaystyle{\prod_{i=0}^{n-1}(-2(n-i))\prod_{i=1}^n(\nu(h)-2(n-i))} & m=2n\\
\displaystyle{\prod_{i=0}^{n-1}(-2(n-i))\prod_{i=0}^n(\nu(h)-2(n-i))} & m=2n+1.
\end{array}
\right.
$$
Then one can easily see that
 \begin{equation}\label{rm}
 \parbox{4.4 in}{\small
  if $w\in M$ is a weight vector of weight $\nu$ with  $xw=0,$ we have
 $x^{m}y^{m}w=r_m(\nu)w.$}
\end{equation}
 We want to show that for each $\lam\in\supp(M),$ $\{k\in\bbbz^{\geq 0}\mid \lam+k\a\in\supp(M)\} $ is bounded. To the contrary, assume there is $\lam\in\supp(M)$ such that  $$\aa:=\{k\in\bbbz^{\geq 0}\mid\lam+k\a\in\supp(M)\}$$ is unbounded. If $\lam(h)$ is not an integer, we set $\mu:=\lam$ and if it is an integer, we pick a  positive integer $m\in \aa$ such that   $(\lam+m\a)(h)$ is positive and set $\mu:=\lam+m\a.$ So in both cases we have
$$\mu(h)+k+2i+1,\mu(h)+k+2i\neq 0 \quad(k\in\bbbz^{>0},\;\; 0\leq i<\frac{k+1}{2}).$$ This implies that
\begin{equation}
\label{non-zero}
r_{k}(\mu+k\a)\neq 0\quad(k\in\bbbz^{>0}).
\end{equation}

 Since $x\in\LL^\a$ acts locally nilpotently and  $\{k\in\bbbz^{>0}\mid \mu+k\a\in\supp(M)\}$ is unbounded, there are $1<k_1<k_2<\cdots$ with $\nu_i:=\mu+k_i\a\in\supp(M)$ and $0\neq v_i\in M^{\nu_i}$ with $xv_i=0.$ Using (\ref{rm}) and (\ref{non-zero}), we get
$$x^{k_i}y^{k_i}v_i=r_{k_i}(\nu_i)v_i\andd r_{k_i}(\nu_i)\neq 0\quad\quad(i\in\bbbz^{>0}).$$
As $y$ acts injectively, $0\neq w_i:=y^{k_i}v_i\in M^\mu.$ But $M^\mu$ is finite dimensional, so one finds $m$ such that $y^{k_{m}}v_{m}=w_{m}=\sum_{i=1}^{m-1}s_iw_{i}=\sum_{i=1}^{m-1}s_iy^{k_{i}}v_{i}$ for some scalars $s_i.$ So we have
\begin{align*}
r_{k_{m}}(\nu_{m})v_{m}=x^{k_{m}}y^{k_{m}}v_{m}=\sum_{i=1}^{m-1}s_ix^{k_{m}}y^{k_{i}}v_{i}=&\sum_{i=1}^{m-1}s_ix^{k_{m}-k_{i}}x^{k_{i}}y^{k_{i}}v_{i}\\
=&\sum_{i=1}^{m-1}r_{k_{i}}(\nu_i)s_ix^{k_{m}-k_{i}}v_{i}=0.
\end{align*}
 But as $r_{k_{m}}(\nu_{m})\neq 0,$ this implies that $v_{m}=0$ which is a contradiction. This completes the proof in the case that $\a\in R^{ln}\cap R_1$. Using the $\frak{sl}_2$-module theory together with the modified argument as above, one can get the result for the  case that $\a\in R^{ln}\cap R_0.$
\qed

\begin{cor}
Suppose that {\rm (s1)} is satisfied for $M,$ then {\rm (s1)} is satisfied for all submodules of $M.$ In particular, if  weight spaces of $M$ are finite dimensional and $M$ has shadow, then  each submodule of $M$  has also  shadow.
\end{cor}
\pf
It is trivial.\qed

\begin{lem}\label{cor1}
Suppose that $M$ has shadow and   $0\neq  \a\in R_{re}.$

\begin{itemize}
\item[(i)]
$\a\in \frak{C}_M$ if and only if $t\a\in {\frak{C}}_M$ for some positive integer $t.$

 \item[(ii)] If either  $\a,-\a\in R^{ln}$ or $\a,-\a\in R^{in},$ then for $\gamma\in R_{re}^\times,$ $\gamma\in R^{in}$ if and only if $r_\a(\gamma)\in R^{in}$ where $r_\a$ is defined as in Lemma \ref{trans}.
     \end{itemize}
\end{lem}
\pf (i) It is trivial using Lemma \ref{trivial} and the fact that  $M$ has shadow.

(ii) If $\a\in R_{re}^\times,$ then   $2\a\in R$ if and only if $\a\in R_1.$ If $\a\in R_1\cap R_{re}^\times,$ then there are  $x\in\LL^\a$ and  $y\in\LL^{-\a}$ such that $$\hbox{span}_\bbbc\{x,y,h:=[x,y],[x,x],[y,y]\}$$ is a Lie superalgebra isomorphic to $\frak{osp}(1,2)$ with $\a(h)=2$ (see \cite[\S~3]{you3} and \cite[Exa. 2.2]{AY}). Then $(\frac{1}{4}[x,x],-\frac{1}{4}[y,y],\frac{1}{2}h)$ is an $\frak{sl}_2$-triple corresponding to $2\a\in R_{re}\cap R_0$ and so  $r_\a=r_{2\a}.$ On the other hand by part (i), $\a\in R^{in}$ if and only if $2\a\in R^{in}.$ So to prove the lemma, without loss of generality, we assume $\a\in R_0.$

We first assume $\pm\a\in R^{ln},$ then we have
\begin{align*}
\gamma\in R^{in} \Longleftrightarrow& \forall \lam\in \hbox{supp}(M)\andd  \forall n\in \bbbz^{\geq0},\; \lam+n\gamma\in \hbox{\rm supp}(M)\\
\stackrel{\tiny {\rm Lem.} \ref{trans}}{\Longleftrightarrow}& \forall \lam\in \hbox{supp}(M)\andd \forall  n\in \bbbz^{\geq0},\; r_\a(\lam)+nr_{\a}(\gamma)\in \hbox{\rm supp}(M)\\
\Longleftrightarrow& r_\a(\gamma)\in R^{in}.
\end{align*}
Next suppose $\pm\a\in R^{in}.$  For  $\gamma\in R_{re}^\times,$ we have $r_\a(\gamma)=\gamma+m\a,$ for some integer $m.$  If $\gamma\in R^{in},$ Lemma \ref{trivial}(ii) implies that $r_\a(\gamma)\in R^{in};$ conversely assume   $r_\a(\gamma)\in R^{in},$ then  by the fact we just proved, $\gamma=r_\a r_\a(\gamma)\in R^{in}.$
\qed
\begin{Thm}\label{close}
Suppose that $M$ is an $\fl$-module having shadow.
Then
\begin{itemize}
\item[(i)] $(R^{ln}+R^{ln})\cap R_{re}^\times\sub R^{ln},$
\item[(ii)] $(R^{ln}+2R^{ln})\cap R_{re}^\times\sub R^{ln}.$
\end{itemize}
\end{Thm}
\pf
(i) Suppose that $\b_1,\b_2\in R^{ln}$ and $\b:=\b_1+\b_2\in R_{re}^\times.$  If $-\b_1\in R^{in},$  then  $\b\in R^{ln}$ as otherwise by Lemma \ref{trivial}(ii), $\b_2=\b-\b_1\in R^{in}$ which is a contradiction. Similarly, if $-\b_2\in R^{in},$ we get $\b\in R^{ln}.$ So to continue the proof, we  assume $\pm\b_1,\pm\b_2\in R^{ln}.$

By Lemma \ref{trivial}, we may assume  $\b_1$ and $\b_2$ are not proportional. Then  either $2(\b_1,\b_2)/(\b_1,\b_1)=\{\pm1,0\}$ or  $2(\b_1,\b_2)/(\b_2,\b_2)=\{\pm1,0\}.$
Without loss of generality, we assume  $2(\b_1,\b_2)/(\b_1,\b_1)=\{\pm1,0\}.$
If $2(\b_1,\b_2)/(\b_1,\b_1)=-1,$ then by Lemma \ref{cor1}(ii), $\b_1+\b_2=r_{\b_1}(\b_2)\in R^{ln}$ and so we are done. So we continue with the case  that $2(\b_1,\b_2)/(\b_1,\b_1)=\{1,0\}.$ Set $r:=2(\b_1,\b_2)/(\b_2,\b_2)$ which is a nonnegative integer. We want to show $\b_1+\b_2\in R^{ln}.$
To the contrary assume  $\b_1+\b_2\in R^{in},$ then by lemma \ref{cor1}(ii), $\b_1-(r+1)\b_2=r_{\b_2}(\b_1+\b_2)\in R^{in}$ and so for each $\lam\in\supp(M)$ and each $k\in\bbbz^{\geq 0},$ using Lemma \ref{trivial}(ii), we have
\begin{align*}
\lam+(r+2)k\b_1=&\lam+k(r+1)(\b_1+\b_2)+k(\b_1-(r+1)\b_2)\\
=&\lam+k(r+1)(\underbrace{\b_1+\b_2}_{\in R^{in}})+k(\underbrace{r_{\b_2}(\b_1+\b_2)}_{\in R^{in}})\in \supp(M)
\end{align*}
 which  contradicts the fact that $\b_1\in R^{ln}\sub \frak{B}_M.$

(ii) Suppose that $\b_1,\b_2,\b_1+2\b_2\in R_{re}^\times$ with $\b_1,\b_2\in R^{ln}.$ If $\b_1+\b_2\in R_{re}^\times,$ we are done using part (i) as $\b_1+2\b_2=(\b_1+\b_2)+\b_2.$ Otherwise,  $\b_1+\b_2\in R_{im}$
and so  $2(\b_1,\b_2)/(\b_2,\b_2)=-2.$ As in part (i), we may assume $\pm\b_2\in R^{ln}.$ Then using Lemma  \ref{cor1}(ii), we have
$\b_1+2\b_2=r_{\b_2}(\b_1)\in R^{ln}.$
\qed

\begin{Thm}\label{property}
Suppose that $M$ is an $\fl$-module having shadow, then for each $\b\in R_{re}^\times,$ one of the following will happen:
\begin{itemize}
\item[\rm (i)] $(\b+\bbbz\d)\cap R\sub R^{ln},$
\item[\rm (ii)] $(\b+\bbbz\d)\cap R\sub R^{in},$
\item[\rm (iii)] there exist $m\in\bbbz$ and $t\in\{0,1,-1\}$ such that for $\gamma:=\b+m\d,$
\begin{align*}
&(\gamma+\bbbz^{\geq 1}\d)\cap R\sub R^{in},\quad (\gamma+\bbbz^{\leq 0}\d)\cap R\sub R^{ln}\\
& (-\gamma+\bbbz^{\geq t}\d)\cap R  \sub R^{in},\quad (-\gamma+\bbbz^{\leq t-1}\d)\cap R\sub R^{ln},
\end{align*}
\item[\rm (iv)]there exist $m\in\bbbz$ and $t\in\{0,1,-1\}$ such that for $\eta:=\b+m\d,$
\begin{align*}
&(\eta+\bbbz^{\leq -1}\d)\cap R\sub R^{in},\quad (\eta+\bbbz^{\geq 0}\d)\cap R\sub R^{ln}\\
& (-\eta+\bbbz^{\leq -t}\d)\cap R  \sub R^{in},\quad (-\eta+\bbbz^{\geq 1-t}\d)\cap R\sub R^{ln}.
\end{align*}
\end{itemize}
\end{Thm}
\pf We know that $\b=\dot\b+n\d$ for some $n\in\bbbz$ and $\dot\b\in \dot R_{re}^\times.$ Using  (\ref{Salpha}), one has $s\in\bbbz^{>0}$ and $k_{\dot\b}\in\bbbz^{\geq0}$ with $\{m\in\bbbz\mid \dot\b+m\d\in R\}=s\bbbz+k_{\dot\b}.$
So $$(\b+\bbbz\d)\cap R=\b+s\bbbz\d.$$
If (i) and (ii) do not hold, then  there is  an integer $k\in \bbbz$ such that
\begin{equation}\label{dag}
\tag{\dag}
\gamma:=\b+sk\d\in R^{ln}\andd \gamma+s\d=\b+sk\d+s\d\in R^{in}
\end{equation} or
\begin{equation}\label{ddag}
\tag{\ddag}
\gamma:=\b+sk\d\in R^{in}\andd \gamma+s\d=\b+sk\d+s\d\in R^{ln}.
\end{equation}
In what follows we show that if (\ref{dag}) (resp. \ref{ddag}) holds, then (iii) (resp. (iv)) is satisfied. We mention that  in (\ref{ddag}), we have $$\eta:=\gamma+s\d\in R^{ln}\andd \eta+s(-\d)=\gamma\in R^{in}.$$ This means that  we just need to study (\ref{dag}). So from now till the end of the proof, we  assume (\ref{dag}) holds. There are four cases:
\begin{itemize}
\item[Case 1.] $-\gamma\in R^{ln}$ and $-\gamma- s\d\in R^{ln}.$
\item[Case 2.] $-\gamma\in R^{in}$ and $-\gamma- s\d\in R^{in}.$
\item[Case 3.] $-\gamma\in R^{in}$ and $-\gamma- s\d\in R^{ln}.$
\item[Case 4.] $-\gamma\in R^{ln}$ and $-\gamma- s\d\in R^{in}.$
\end{itemize}

\noindent {\bf Case 1.} In this case, we have  $\pm\gamma\in R^{ln}.$ So Lemma  \ref{cor1} implies that
\begin{equation}
\label{symm}
\gamma+p s\d\in R^{in}\Leftrightarrow -\gamma+p s\d\in R^{in}\quad\quad (p\in\bbbz).
\end{equation}
In particular, since ($\dag$) holds, we have  $\gamma+ s\d\in R^{in}$ and so  $-\gamma+ s\d\in R^{in}.$
In two steps we show the following:
{\small
\begin{align*}\label{shar1}
\tag{$\sharp_1$}
  \overbrace{\cdots\quad  \gamma-2s\d\quad  \gamma-s\d \quad\gamma}^{\in R^{ln}} \quad\overbrace{\gamma+s\d \quad \gamma+2s\d\quad \cdots}^{\in R^{in}}
  \\
  \underbrace{\cdots\quad  -\gamma-2s\d\quad  -\gamma-s\d \quad-\gamma}_{\in R^{ln}} \quad\underbrace{-\gamma+s\d \quad -\gamma+2s\d \quad\cdots}_{\in R^{in}}
\end{align*}
}

\noindent\underline{Claim 1.} For $n\in\bbbz^{\geq 1},$ we have $\pm\gamma+ns\d\in R^{in}:$ Let $n\in\bbbz^{\geq 1},$ then  by Lemma \ref{trivial}(ii)
{\small $$\pm\gamma+ (1+2n) s\d=(\underbrace{\pm\gamma+ s\d}_{\in R^{in}})+n(\underbrace{\gamma+ s\d}_{\in R^{in}})+n(\underbrace{-\gamma+ s\d}_{\in R^{in}})\in R^{in}.$$} Also we have
{\small \begin{align*}
\pm \gamma+2n s\d=&(\pm\gamma+2 s\d)+(n-1)(\underbrace{\gamma+ s\d}_{\in R^{in}})+(n-1)(\underbrace{-\gamma+ s\d}_{\in R^{in}})
\end{align*}}
 which is an element of $R^{in}$ provided that $\pm\gamma+2 s\d\in R^{in}.$ If  to the contrary $\pm\gamma+2 s\d\in R^{ln},$ then by Theorem \ref{close}(ii) $$\gamma+3 s\d=(\underbrace{-\gamma- s\d}_{\in R^{ln}})+2(\underbrace{\gamma+2 s\d}_{\in R^{ln}})\in R^{ln}$$ while   {\small $$-\gamma+3 s\d=(\underbrace{-\gamma+ s\d}_{\in R^{in}})+(\underbrace{\gamma+ s\d}_{\in R^{in}})+(\underbrace{-\gamma+ s\d}_{\in R^{in}})\in R^{in}$$} which  contradicts (\ref{symm}). This completes the proof in of Claim 1.

\noindent\underline{Claim 2.} For all positive integers $n,$  $\pm\gamma-ns\d\in R^{ln}:$ If $n$ is a positive integer with $\pm\gamma-2n s\d\in R^{in},$ then {\small $$\pm\gamma=(\underbrace{\pm\gamma-2n s\d}_{\in R^{in}})+n(\underbrace{\gamma+ s\d}_{\in R^{in}})+n(\underbrace{-\gamma+ s\d}_{\in R^{in}})\in R^{in}$$} which is a contradiction. Also if  $\pm\gamma+(-2n-1) s\d\in R^{in}$ for some nonnegative integer $n,$  then {\small $$\pm\gamma- s\d=(\underbrace{\pm\gamma+(-2n-1) s\d}_{\in R^{in}})+n(\underbrace{\gamma+ s\d}_{\in R^{in}})+n(\underbrace{-\gamma+ s\d}_{\in R^{in}})\in R^{in}$$} which  contradicts our assumption in Case 1; see (\ref{symm}). This completes the proof.

\noindent {\bf Case 2.} In this case we show:
{\small \begin{align*}\label{shar2}
\tag{$\sharp_2$}
  \overbrace{\cdots\quad  \gamma-2s\d\quad  \gamma-s\d \quad\gamma}^{\in R^{ln}} \quad\overbrace{\gamma+s\d \quad \gamma+2s\d\quad \cdots}^{\in R^{in}}
  \\
  \underbrace{\cdots\quad  -\gamma-2s\d}_{\in R^{ln}}\quad  \underbrace{-\gamma-s\d \quad-\gamma \quad-\gamma+s\d \quad -\gamma+2s\d\quad \cdots}_{\in R^{in}}
\end{align*}}

\noindent\underline{Claim 1.} For all nonnegative integers $n,$ $\gamma-n s\d\in R^{ln}$:  Suppose to the contrary that $n$ is a positive integer and $\gamma-n s\d\in R^{in}, $ using (\ref{dag}), we have {\small $$\underbrace{\gamma}_{\in R^{ln}}=(\underbrace{\gamma-n s\d}_{\in R^{in}})+n(\underbrace{-\gamma}_{\in R^{in}})+n(\underbrace{\gamma+ s\d}_{\in R^{in}})\in R^{in}$$} which is a contradiction.

\noindent\underline{Claim 2.} For  $n\in\bbbz^{\geq2},$ $-\gamma-n s\d\in R^{ln}$: We first note that as $\pm(\gamma+ s\d)\in R^{in}$ (by (\ref{dag}) and our assumption), then by Lemma \ref{cor1}, $-\gamma-2 s\d=\gamma-2\gamma-2 s\d=r_{\gamma+ s\d}(\gamma)\in R^{ln}.$ Now if
 to the contrary, for some  $n\in \bbbz^{\geq3},$  $-\gamma-n s\d\in R^{in},$ then {\small $$\underbrace{-\gamma-2 s\d}_{\in R^{ln}}=(\underbrace{-\gamma-n s\d}_{\in R^{in}})+(n-2)(\underbrace{-\gamma}_{\in R^{in}})+(n-2)(\underbrace{\gamma+ s\d}_{\in R^{in}})\in R^{in}$$} which is a contradiction.

\noindent\underline{Claim 3.} For all $n\in \bbbz^{\geq -1},$ we have $\gamma+(n+2) s\d,-\gamma+n s\d\in R^{in}$: By our assumption in Case 2 and (\ref{dag}), $-\gamma,-\gamma-s\d,\gamma+s\d\in R.$ Also if $n$ is  a nonnegative  integer, then
{\small \begin{align*}
&-\gamma+n s\d=(n+1)(\underbrace{-\gamma}_{\in R^{in}})+n(\underbrace{\gamma+ s\d}_{\in R^{in}})\in R^{in},\\
&\gamma+(n+2) s\d=(n+1)(\underbrace{-\gamma}_{\in R^{in}})+(n+2)(\underbrace{\gamma+ s\d}_{\in R^{in}})\in R^{in}.
\end{align*}}
\noindent{\bf Case 3.}  We shall show the following:
\begin{align*}\label{shar3}
\tag{$\sharp_3$}
  \overbrace{\cdots\quad  \gamma-2s\d\quad  \gamma-s\d \quad\gamma}^{\in R^{ln}} \quad\overbrace{\gamma+s\d \quad \gamma+2s\d \quad \cdots}^{\in R^{in}}
  \\
  \underbrace{\cdots\quad  -\gamma-2s\d\quad  -\gamma-s\d }_{\in R^{ln}}\quad\underbrace{-\gamma \quad-\gamma+s\d \quad -\gamma+2s\d\quad \cdots}_{\in R^{in}}
\end{align*}

\noindent\underline{Claim 1.} For all nonnegative integers $n,$ $-\gamma+n s\d,\gamma+(n+1)s\d\in R^{in}$: Suppose that $n\geq 0,$  then
\begin{align*}
&-\gamma+n s\d=(n+1)(\underbrace{-\gamma}_{\in R^{in}})+n(\underbrace{\gamma+ s\d}_{\in R^{in}})\in R^{in},\\
&\gamma+(n+1) s\d=n(\underbrace{-\gamma}_{\in R^{in}})+(n+1)(\underbrace{\gamma+ s\d}_{\in R^{in}})\in R^{in}.
\end{align*}
This completes the proof.

\noindent\underline{Claim 2.}
For all nonnegative integers $n,$ $\gamma-n s\d\in R^{ln}$:  We know from (\ref{dag}) that $\gamma\in R^{ln}.$ Suppose to the contrary that $n$ is a positive integer and $\gamma-n s\d\in R^{in}.$ As by Claim 1, $-(\gamma-n s\d)\in R^{in},$  we have using Lemma \ref{cor1} that
$$-\gamma+2n s\d=\gamma-2\gamma+2n s\d=r_{\gamma-n s\d}(\gamma)\in R^{ln}$$ which   contradictions Claim 1.

\noindent\underline{Claim 3.} For all positive integers $n,$ $-\gamma-n s\d\in R^{ln}$: By our assumption, $-\gamma- s\d\in R^{ln}.$ So using Claim 2 and Lemma \ref{close}, we have
\begin{align*}
-\gamma-(n+1) s\d=(\gamma-(n-1) s\d)+2(-\gamma- s\d)\in R^{ln}+2R^{ln}\sub R^{ln}.
\end{align*}

\noindent{\bf Case 4.}
We show that this case cannot happen.  If $-\gamma\in R^{ln}$ and $-\gamma-s\d\in R^{in},$ by (\ref{dag}), we have  $\pm\gamma\in R^{ln}$ and $\pm(\gamma+ s\d)\in R^{in}.$ So Lemma \ref{cor1} implies that  $\pm\gamma\pm s\d=r_{\gamma}(\pm (\gamma+ s\d))\in R^{in}.$ In particular
\begin{equation}\label{if-and-if}
\parbox{4.5in}{$\mu +(-\gamma+ s\d)\in\supp(M)\Leftrightarrow\mu\in \supp(M)\Leftrightarrow \mu +(\gamma+ s\d)\in\supp(M).$ }
\end{equation}

Now suppose $\lam\in \supp(M).$ Since $\gamma\in R^{ln},$ we find a positive integer $p$ such that $\lam+2p\gamma\not\in\supp(M).$ So
{\small \begin{align*}
\lam+2p\gamma\not\in\supp(M)\stackrel{(\ref{if-and-if})}{\Longrightarrow}&\lam+2p\gamma+2p(-\gamma+ s\d)\not\in \supp(M)\\
\Longrightarrow&\lam+2p s\d\not\in \supp(M)\\
\Longrightarrow&\lam+p(\gamma+ s\d)+p(-\gamma+ s\d)\not\in\supp(M)
\stackrel{(\ref{if-and-if})}{\Longrightarrow}\lam\not\in \supp(M).
\end{align*}}
This is a contradiction.\qed

\begin{deft}\label{full-hybrid}
{\rm Suppose that $M$ is an $\fl$-module having shadow. We say $\a\in R_{re}^\times$ is {\it full-locally nilpotent} (resp. {\it full-injective}) if $(\a+\bbbz\d)\cap R\sub R^{ln}$  (resp. $(\a+\bbbz\d)\cap R\sub R^{in}$), otherwise, we call it {\it hybrid}.
}
\end{deft}

\section{Modules over twisted affine Lie superalgebras}

Keeping the same notations as in previous sections, throughout this section, we assume  $\fl$ is a twisted affine Lie superalgebra of type $X=A(2k-1,2\ell-1)^{(2)}$ ({\tiny$(k,\ell)\neq (1,1)$}), $A(2k,2\ell)^{(4)},$    $A(2k,2\ell-1)^{(2)}$ and $D(k,\ell)^{(2)}$ where  $k,\ell$ are positive integers, $\fh\sub \fl_0$ is the standard Cartan subalgebra of $\fl$ with corresponding root system $R;$ see Table~\ref{table1}, and   $R_0$ (resp. $R_1$)   is the set of weights of $\fl_0$ (resp. $\fl_1$) with respect to $\fh.$

 Assume $M$ is an irreducible  $\fl$-module having a weight space decomposition with respect to $\fh$ with finite dimensional weight spaces. By Proposition \ref{suff}, $M$ has shadow.  We know from (\ref{Salpha}) that there is  $r\in\bbbz^{>0}$ such that
 \begin{equation}
 \label{r}
 R_i+r\bbbz\d\sub R_i\quad\quad(i=0,1).
 \end{equation}
Since $M$ has shadow, using  Theorem \ref{property}, we have
$$\hbox{\small $R_{re}^\times= \underbrace{\{\a\in R_{re}^\times\mid \exists N, (\a+\bbbz^{\geq N}\d)\cap R\sub  R^{ln}\}}_{K_1}\uplus  \underbrace{\{\a\in R_{re}^\times\mid \exists N, (\a+\bbbz^{\geq N}\d)\cap R\sub  R^{in}\}}_{K_2}$ }$$ in which ``~$\uplus$~'' indicates disjoint union.
 If $\a,\b\in K_1$ (resp. $\in K_2$) and $\a+\b\in R_{re}^\times, $ then for large enough $n,$ (\ref{r}) implies that  $\a+nr\d,\b+rn\d\in R^{ln}$ (resp. $\in  R^{in}$) and by Theorem \ref{close} (resp. Lemma~\ref{trivial}(ii)), $\a+\b+2rn\d\in R^{ln}$ (resp. $\a+\b+2rn\d\in R^{in}$); i.e., $\a+\b\in K_1$ (resp. $\in K_2$). It means that
 \begin{equation}\label{closed}
 \parbox{2.3in}{$K_1$ and $K_2$ are closed subsets of $R_{re}^\times.$}
\end{equation}
We know from remark \ref{decom-aff-app} that there are  affine Lie subalgebras $\fl_0(1)$ and $\fl_0(2)$ of $\fl_0$ with Cartan subalgebras $\fh_1$ and $\fh_2$ respectively such that $$\fh=\fh_1+\fh_2 .$$
Set
\begin{equation*}\label{ki}
\fk_i:=\fl_0(i)+ \fh\quad\quad(i=1,2).
\end{equation*}
We denote by $R(i),$  the set of weights of $\fk_i$ with respect to $\fh;$ this is in fact the root system of $\fl_0(i)$ with respect to $\fh_i.$

\begin{lem}
\label{new5}
 Suppose that $i\in\{1,2\}$ and $R(i)=R(i)^+\cup R(i)^\circ\cup R(i)^-$ is a triangular decomposition for $R(i)$ with corresponding functional $\boldsymbol\zeta$ such that $R(i)^+\cap R_{re}\sub R^{ln}$ and $R(i)^-\cap R_{re}\sub R^{in}.$  Assume  $\boldsymbol\zeta(\d)>0$ and $W$ is an $\fl_0$-submodule of $M,$ then there is a positive integer $p$  and $\lam\in \supp(W)$ with $(\lam+\bbbz^{>0}p\d)\cap\supp(W)=\emptyset.$
\end{lem}
\pf  It essentially follows from \cite[\S 2]{DG} but for the convenience of readers, we give the proof. Since $\boldsymbol\zeta(\d)>0,$ it follows that  $R(i)^\circ$ is either $\{0\}$ or a finite root system and that if $\Sigma$ is the standard base of $R(i),$ the set of positive roots of $R(i)$ with respect to $\Sigma$ intersects $R(i)^-$ in a   finite set.  So, by \cite[Pro. 2.10(i)]{DFG}, there is a base  $B=\{\a_1,\ldots,\a_\ell\}$  of $R(i)$ contained in $P:=R(i)^+\cup R(i)^\circ.$
 If $R(i)^\circ\neq \{0\},$ set
 $$B_1:=B\cap R(i)^+$$ and assume $\w$  is the Weyl group of the finite root system $R(i)^\circ.$    We set
 $\Phi:=\w(B_1)\sub R(i)^+\cap R_{re}.$  Then there is $p\in \bbbz^{\gneq0}$ such that $p\d\in \sspan_{\bbbz^{\geq 0}}\Phi;$ see (2.15) of \cite{DG}. Moreover, using \cite[Pro.~2.1.1]{MP} and  Lemma \ref{general}(ii) together with the fact that $M$ has shadow, there is $\lam\in \supp(W)$ such that $(\lam+\hbox{span}_{\bbbz^{\geq 0}}\Phi)\cap \supp(W)=\{\lam\}.$ So $(\lam+\hbox{span}_{\bbbz^{\geq 0}}p\bbbz\d)\cap \supp(W)=\{\lam\}$ as we desired.

    Next assume $R(i)^\circ=\{0\}.$ Therefore,   we have $B\sub R(i)^+\cap R_{re}$ and so by Lemma \ref{general}(ii), there is $\lam\in \supp(W)$ such that $(\lam+\hbox{span}_{\bbbz^{\geq 0}}B)\cap \supp(W)=\{\lam\}.$ But $R(i)_{im}=s\bbbz\d$ for some positive integer $s$ and as $\boldsymbol\zeta(\d)>0$ and $B\sub R(i)^+,$ we have  $s\d\in \hbox{span}_{\bbbz^{\geq 0}}B.$ This completes the proof.\qed

\bigskip

We set $$R(i)^{ln}:=R(i)\cap R^{ln}\andd R(i)^{in}:=R(i)\cap R^{in}.$$
\begin{deft}
{\rm
\begin{itemize}
\item[(i)]
 We say $R(i)$ is {\it tight} if there is a  nonzero real root $\a\in R(i)$ with $(\a+\bbbz\d)\cap R(i)\sub R(i)^{ln}$ or   $(\a+\bbbz\d)\cap R(i)\sub  R(i)^{in};$ otherwise, we call it {\it hybrid}.
\item[(ii)] We say $M$ is {\it hybrid} if both $R(1)$ and $R(2)$ are hybrid; otherwise, we call it {\it tight.}
\end{itemize}
  }
\end{deft}

If $R(i)$ is hybrid, (\ref{closed})  together with Theorem \ref{property} implies that $R(i)\cap K_1$ as well as $R(i)\cap K_2$ are symmetric closed subsets of $R(i)_{re}^\times$ which in turn implies that
$(\a,\b)=0$ if $\a\in R(i)\cap K_1$ and $\b\in R(i)\cap K_2.$
Therefore, either $R(i)\cap K_1=\emptyset$ or $R(i)\cap K_2=\emptyset$ as $R(i)$ is an affine root system.

\begin{deft}
{\rm Suppose that  $R(i)$ ($i=1,2$) is hybrid. We call $R(i)$   {\it up-nilpotent hybrid} if  $R(i)\cap K_1=R(i)^\times_{re},$ otherwise, we call it   {\it down-nilpotent hybrid}. We set {\small
\begin{equation}\label{parabolic}
P_i:=\left\{
\begin{array}{ll}
R(i)^{ln}\cup -R(i)^{in}\cup (\bbbz^{\geq 0}\d\cap R(i))&\hbox{if $R(i)$ is up-nilpotent hybrid }\\
R(i)^{ln}\cup -R(i)^{in}\cup (\bbbz^{\leq 0}\d\cap R(i)) &\hbox{if $R(i)$ is down-nilpotent hybrid}.
\end{array}
\right.
\end{equation}
}
}
\end{deft}

\begin{lem}\label{para}
Suppose that  $R(i)$ ($i=1,2$) is hybrid. Then $P_i$ is  a proper parabolic subset of  $R(i);$ i.e., $P_i$ is a proper subset of $R(i)$ satisfying $R(i)=P_i\cup-P_i$ and $(P_i+P_i)\cap R(i)\sub P_i.$
\end{lem}
\pf It is trivial that $P_i$ is proper. Also as $R(i)=P_i\cup -P_i,$  we just need to show that $P_i$ is closed. We first assume $R(i)$ is down-nilpotent hybrid.  Using Theorem \ref{close}, Lemma \ref{trivial}(ii) as well as  Theorem \ref{property} and (\ref{shar1})-(\ref{shar3}) in its proof, we get
$$((R(i)^{ln}\cup -R(i)^{in})+(R(i)^{ln}\cup -R(i)^{in}))\cap R(i)\sub P_i.$$ So we just need to prove  $R(i)\cap(( R(i)^{ln}\cup -R(i)^{in})+(\bbbz^{<0}\d\cap  R(i)))\sub P_i.$ Suppose $\a\in -R(i)^{in}$  and $m\in\bbbz^{<0}$ are such that $\a+m\d\in R(i),$ then as  $\a\in -R(i)^{in},$ Theorem \ref{property} implies that $-\a-m\d\in R(i)^{in}$ and so $\a+m\d\in -R(i)^{in}.$ Similarly, we can see that     $\a+m\d\in R(i)^{ln}$ if $\a\in R(i)^{ln}$  and $m\in\bbbz^{<0}$ with $\a+m\d\in R(i).$  Using the same argument as above, one can get the result when $R(i)$ is up-nilpotent hybrid.\qed

\begin{rem}\label{rem2}
{\rm Suppose $i=1,2$ and   $s_i$ is the positive integer  with
\[R(i)_{im}=s_i\bbbz\d.\] Assume $R(i)$ is up-nilpotent hybrid, so  we have
\begin{equation}\label{si}
s_i\d\in P_i\setminus-P_i.
\end{equation} One knows  from affine Lie theory that each base of $R(i) $ is of the form  $\pm\Sigma_i$ for  $$\Sigma_i:=\{\dot\b_1,\ldots,\dot\b_t,s_i\d-\theta_i\}$$ where    $B_i:=\{\dot\b_1,\ldots,\dot\b_t\}$  is a base of an irreducible  finite root system $\dot R_i$ with
\begin{equation}\label{new7}
(\dot R_i)_{ind}\cup((\dot R_i)_{sh}+s_i\bbbz\d)\sub R(i),
\end{equation} in which
\begin{center}
\parbox{3.4in}{\small $(\dot R_i)_{sh}$  is the set of short roots, that is, the set of roots of $\dot R_i$ with  the smallest length and $(\dot R_i)_{ind}=(\dot R_i\setminus 2\dot R_i)\cup\{0\}$,}
\end{center}
 and $\theta_i$ is as in the following table:

\begin{table}[h]\caption{\small Description of $\theta_i$}\label{highest}
\begin{tabular}{|c|c|}
\hline
\small{Type of $R(i)$} & $\theta_i$\\
\hline
$\hbox{untwisted types}$  & \small{The highest root of $\dot R_i$ with respect to $B_i$}\\
\hline
 $A_{2p}^{(2)}$ $(p\geq 1)$& \small{2 times of the highest short root of $\dot R_i$ with respect to $B_i$}\\
\hline
$\hbox{other types}$& \small{The highest short root of $\dot R_i$ with respect to $B_i$}\\
\hline
\end{tabular}
\end{table}

Here, we use affine's labels from Kac's Book \cite{K2}.
In particular,
\begin{equation}\label{new6}
\parbox{4.3in}{\small $\frac{1}{2}\theta_i\in (\dot R_i)_{ind}\sub R(i)$ if $R(i)$ is  of type $A_{2p}^{(2)}$ and $\theta_i\in (\dot R_i)_{ind}\sub R(i),$ otherwise. }
\end{equation}
 Moreover, each positive root of $R(i)$ with respect to $\Sigma_i$ either is a positive root of $(\dot R_i)_{ind}$ with respect to $B_i$ or is of the form $\dot\a+m\d,$ for some  root $\dot\a\in \dot R$  and a positive integer $m.$ This together with Proposition~2.10 of \cite{DFG}, (\ref{si}) and  the fact that for each $\a\in R(i),$ $(\a+\bbbz\d)\cap P_i\neq \emptyset,$
 implies that
 \begin{equation}\label{new8}
 \parbox{4.4in}{
 $\bullet$ there is a base $\Pi_i$ of $R(i)$ such that the set $R(i)^+(\Pi_i)$ of positive roots of $R(i)$ with respect to $\Pi_i$ is a subset of $ P_i.$\\
 $\bullet$ there is a functional $\boldsymbol\zeta$ on $\sspan_\bbbr R(i)$ with  $P_i=R(i)^+\cup R(i)^\circ$ and $\boldsymbol\zeta(\d)>0;$ see Definition~\ref{def1}.
 }
 \end{equation}

  We claim that  $\Pi_i$ is of the form $\Sigma_i.$ To the contrary, assume $\Pi_i$ is of the form $-\Sigma_i.$ So there is a finite root system $\dot R_i$ satisfying (\ref{new7})
 and a  base $\{\dot\b_1,\ldots,\dot\b_t\}$ of $\dot R_i$ such that
$$\Pi_i=\{-\dot\b_1,\ldots,-\dot\b_t,\theta_i-s_i\d\}\sub P_i$$ where $\theta_i$ is as in Table~\ref{highest}.

We first assume $R(i)$ is of type $A_{2p}^{(2)}$ $(p\geq 1)$. Contemplating (\ref{new6}), as  $\dot R_i$ is a finite root system and $\{-\dot\b_1,\ldots,-\dot\b_t\}\sub P_i,$ we get  $-\dot\b:=-\frac{1}{2}\theta_i\in P_i.$
 Also we know from (\ref{new7})  and Table~\ref{highest} that    $\dot\b-s_i\d\in R(i).$
So we have
  \begin{align*}
 \hbox{\small$ -s_i\d=(2\dot\b-s_i\d)+(-\dot\b)+(-\dot\b)=\underbrace{(\theta_i-s_i\d)+(-\dot\b)}_{\in (P_i+P_i)\cap R(i)}+(-\dot\b)\in (P_i+P_i)\cap R(i)\sub P_i,$}
  \end{align*}
in other words, $s_i\d\in P_i\cap -P_i$ which is a contradiction.

Also if  $R(i)$ is not of type $A_{2p}^{(2)},$ then using Table~\ref{highest}, we have $-\theta_i\in P_i$ and so
$$-s_i\d=(\theta_i-s_i\d)+(-\theta_i)\in (P_i+P_i)\cap R(i)\sub P_i$$ which is again a contradiction.
}
\end{rem}
\begin{lem}\label{same-type}
Suppose that $j,{j'}\in\{1,2\}$ and $j\neq {j'}.$ If $R(j)$ is up-nilpotent hybrid  (resp. down-nilpotent hybrid), then $R({j'})$ is  either tight or up-nilpotent hybrid (resp. down-nilpotent hybrid).
\end{lem}

\pf To the contrary, assume $R(j)$ is up-nilpotent hybrid and $R({j'})$
is down-nilpotent hybrid.  By (\ref{new8}), there is a functional $\boldsymbol\zeta$ on $\sspan_\bbbr R(j)$ with  $P_j=R(j)^+\cup R(j)^\circ$ and $\boldsymbol\zeta(\d)>0.$ Using  Lemma~\ref{new5}, one  finds $p\in\bbbz^{>0}$ and $\mu\in \supp(M)$ such that
\begin{equation}\label{p}
(\mu +\bbbz^{>0}p\d)\cap \supp(M)=\emptyset.
\end{equation} For $r$ as in (\ref{r}) and $\b\in R({j'})_{re}^\times,$ since $R({j'})$ is down-nilpotent hybrid, we pick $m>0$ such that
$$\pm\b-nrp\d\in R({j'})^{ln}\andd \pm\b+nrp\d\in R({j'})^{in}\quad(n\geq m).$$
Now if $\mu +\b-mrp\d\in \supp(M),$ then as $-\b+2mrp\d\in R^{in},$ we have $$\mu+mrp\d=(\mu+\b-mrp\d)-\b+2mrp\d\in \supp(M)$$ which is a contradiction due to (\ref{p}), in particular, $$(\fk_{j'})^{\b-mrp\d}M^\mu=\{0\}.$$
Also as  $\b,\b+2mrp\d\in R({j'})_{re}^\times,$ the root string property for the affine root system $R(j')$ implies that $2mrp\d\in R({j'})$ and by (\ref{p}), we have   $$(\fk_{j'})^{2mrp\d}M^\mu=\{0\}.$$
Therefore, we have $$(\fk_{j'})^{\b+mrp\d} M^\mu=[(\fk_{j'})^{\b-mrp\d},(\fk_{j'})^{2mrp\d}]M^\mu=\{0\}$$ which  contradictions the fact that $\b+mrp\d\in R({j'})^{in}$.\qed

\begin{lem}\label{bothhybrid}
Suppose that $R(1)$ and $R(2)$ are hybrid and recall  (\ref{parabolic}). Set $P:=P_1\cup P_2.$ Then there exists a functional $\boldsymbol\zeta:\sspan_\bbbr R_0\longrightarrow \bbbr$ such that $$P=\{\a\in R_0\mid \boldsymbol\zeta(\a)\geq 0\};$$
in particular,
\[\{\a\in R_0\cap R_{re}\mid \boldsymbol\zeta(\a)>0\}\sub R^{ln}\andd \{\a\in R_0\cap R_{re}\mid \boldsymbol\zeta(\a)<0\}\sub R^{in}.\]
\end{lem}
\pf
Without loss of generality, using Lemma  \ref{same-type}, we assume both $R(1)$ and $R(2)$ are up-nilpotent hybrid.
We use  Remark~\ref{rem2} to choose  bases $\Pi_1$ and $\Pi_2$ of  $R(1)$ and $R(2)$ respectively as
$$\hbox{\small $\Pi_1=\{\a_j,\a_0:=s_1\d-\theta_1\mid 1\leq j\leq n\}\sub P_1,\;\;\Pi_2=\{\b_j,\b_0:=s_2\d-\theta_2\mid 1\leq j\leq m\}\sub P_2$}$$ in which $s_1$ and $s_2$ are defined by
\[R(i)_{im}=s_i\bbbz\d\quad\quad(i=1,2),\]
$B_1:=\{\a_1,\ldots,\a_n\}$ and  $B_2:=\{\b_1,\ldots,\b_m\}$ are  bases of some  finite root systems $\dot R_1$ and $\dot R_2$ with  $(\dot R_1)_{ind}\sub R(1)$ and $(\dot R_2)_{ind}\sub R(2)$ respectively  and  $\theta_i $ ($i=1,2$) is as in Table~\ref{highest}.
Renumbering the elements of $ B_1$ and $B_2$ if necessary,  we assume $$\a_1,\ldots,\a_t,\b_1,\ldots,\b_k\in P\setminus -P\andd \a_{t+1},\ldots,\a_n,\b_{k+1},\ldots,\b_m\in P\cap -P.$$ Using a modified argument as in \cite[Pro. 2.10(ii)]{DFG}, we just need to define a functional $\boldsymbol\zeta$ satisfying
\begin{equation}
\label{main}
\boldsymbol\zeta(\Pi_i\cap (P\setminus-P))\sub \bbbr^{>0}\andd \boldsymbol\zeta(\Pi_i\cap (P\cap-P))=\{0\}\;\;\;\;(i=1,2).
\end{equation}
Since $B:=\Pi_1\cup\Pi_2\setminus\{s_2\d-\theta_2\}$ is a basis for the vector space  $\sspan_{\bbbr}R_0,$ to define $\boldsymbol\zeta,$ it is enough to define $\boldsymbol\zeta$ on  $B.$
Let
$$\theta_1=\sum_{i=1}^nr_i\a_i\andd \theta_2=\sum_{j=1}^mk_j\b_j$$ and recall from finite dimensional Lie theory
 that $r_i$'s as well as $k_i$'s are positive integers. We then set $$s:=s_2/s_1.$$

\noindent {\bf Case 1. $\pmb{s_1\d-\theta_1,s_2\d-\theta_2\in P\cap -P:}$}  Define
\begin{align*}
\boldsymbol\zeta:& \sspan_{\bbbr}R_0\longrightarrow \bbbr;\;\;\;\;\;\left\{
\begin{array}{rll}
s_1\d-\theta_1&\mapsto 0 & \\
\a_i&\mapsto \frac{1}{str_i}& 1\leq i\leq t,\\
\a_i&\mapsto 0& t+1\leq i\leq n,\\
\b_j&\mapsto \frac{1}{kk_j}& 1\leq j\leq k,\\
\b_j&\mapsto 0& k+1\leq j\leq m.\\
\end{array}
\right.
\end{align*}
Then
$$\boldsymbol\zeta(s_2\d-\theta_2)=s\boldsymbol\zeta(s_1\d)-\boldsymbol\zeta(\theta_2)=
s\boldsymbol\zeta(s_1\d-\theta_1)+s\boldsymbol\zeta(\theta_1)-\boldsymbol\zeta(\theta_2)=s\boldsymbol\zeta(s_1\d-\theta_1)=0.$$

\smallskip

\noindent {\bf Case 2. $\pmb{s_1\d-\theta_1,s_2\d-\theta_2\in P\setminus -P:}$}
Define
\begin{align*}
\boldsymbol\zeta:& \sspan_{\bbbr}R_0\longrightarrow \bbbr;\;\;\;\;\;
 \left\{
\begin{array}{rll}
s_1\d-\theta_1&\mapsto 1 & \\
\a_i&\mapsto \frac{1}{str_i}& 1\leq i\leq t,\\
\a_i&\mapsto 0& t+1\leq i\leq n,\\
\b_j&\mapsto \frac{1}{kk_j}& 1\leq j\leq k,\\
\b_j&\mapsto 0& k+1\leq j\leq m.\\
\end{array}
\right.
\end{align*}
Then $$\boldsymbol\zeta(s_2\d-\theta_2)=s\boldsymbol\zeta(s_1\d)-\boldsymbol\zeta(\theta_2)=s\boldsymbol\zeta(s_1\d-\theta_1)+s\boldsymbol\zeta(\theta_1)-\boldsymbol\zeta(\theta_2)=s\boldsymbol\zeta(s_1\d-\theta_1)=s.$$

\noindent {\bf Case 3. $\pmb{s_1\d-\theta_1\in P\setminus -P\andd s_2\d-\theta_2\in P\cap -P:}$} Define \begin{align*}
\boldsymbol\zeta:& \sspan_{\bbbr}R_0\longrightarrow \bbbr;\;\;\;\;\;\left\{
\begin{array}{rll}
s_1\d-\theta_1&\mapsto {\frac{1}{s}} & \\
\a_i&\mapsto \frac{1}{str_i}& 1\leq i\leq t,\\
\a_i&\mapsto 0& t+1\leq i\leq n,\\
\b_j&\mapsto \frac{2}{kk_j}& 1\leq j\leq k,\\
\b_j&\mapsto 0& k+1\leq j\leq m.\\
\end{array}
\right.
\end{align*}
Then
$$\boldsymbol\zeta(s_2\d-\theta_2)=s\boldsymbol\zeta(s_1\d-\theta_1)+s\boldsymbol\zeta(\theta_1)-\boldsymbol\zeta(\theta_2)=1+1-2=0.$$

\noindent {\bf Case 4. $\pmb{s_1\d-\theta_1\in P\cap -P\andd s_2\d-\theta_2\in P\setminus -P:}$} Define \begin{align*}
\boldsymbol\zeta:& \sspan_{\bbbr}R_0\longrightarrow \bbbr;\;\;\;\;\;\left\{
\begin{array}{rll}
s_1\d-\theta_1&\mapsto 0 & \\
\a_i&\mapsto \frac{1}{str_i}& 1\leq i\leq t,\\
\a_i&\mapsto 0& t+1\leq i\leq n,\\
\b_j&\mapsto \frac{1}{2kk_j}& 1\leq j\leq k,\\
\b_j&\mapsto 0& k+1\leq j\leq m.\\
\end{array}
\right.
\end{align*}
Then
$$\boldsymbol\zeta(s_2\d-\theta_2)=s\boldsymbol\zeta(s_1\d-\theta_1)+s\boldsymbol\zeta(\theta_1)-\boldsymbol\zeta(\theta_2)=0+1-\frac{1}{2}=\frac{1}{2}.$$

This completes the proof.
\qed

\begin{Thm}\label{triangular}
Suppose that $R(1)$ and $R(2)$ are hybrid. Then there is a triangular decomposition $R=R^+\cup R^\circ\cup R^-$ for $R$ such that $$M^{\fl^+}=\{v\in M\mid \fl^\a v=\{0\}\;\; (\forall \a\in R^+)\}\neq \{0\}.$$
\end{Thm}
\pf Without loss of generality, we assume both $R(1)$ and $R(2)$ are up-nilpotent hybrid and define the functional $\boldsymbol\zeta:\sspan_\bbbr R_0\longrightarrow \bbbr$ as in Lemma \ref{bothhybrid}. Since  $\sspan_\bbbr R_0=\sspan_\bbbr R$ (see (\ref{span})), $\boldsymbol\zeta$ defines  a triangular decomposition   $R=R^+\cup R^\circ\cup R^-$ for $R.$ We note that as two times of a real odd root is a real even root, Lemma~\ref{trivial} and Theorem~\ref{close} imply that  $$\d\in R^+,\;\; R^+\cap R_{re}\sub R^{ln}\andd  R^-\cap R_{re}\sub R^{in}.$$
We set
\begin{align*}
\aa:=&\{v\in M\setminus\{0\}\mid \fl^{\a}v= \{0\} \;\;\forall\a\in  R^+\cap( R_{re}\cup R_{im})\}\\
=&\{v\in M\setminus\{0\}\mid \fl^{n\d}v= \fl^{\a}v= \{0\} \;\;\forall\a\in  R_{re}\cap R^+,\; n\in\bbbz^{>0}\}.
\end{align*}
 Then using Proposition \ref{last-one},  it is enough to show
\begin{equation}\label{what-we-need}
\parbox{3.4in}{there exists $v\in \aa$ such that for each $\dot\a\in \dot R_{ns}^\times,$ there is $N\in\bbbz^{\geq0}$ with $\fl^{\dot \a+n\d} v=\{0\}$  for all $n\geq N.$  }
\end{equation}
 Apply Lemma~\ref{new5} to find a positive integer $p$ and $\lam\in \supp(M)$ such that $(\lam+{\bbbz^{>0}}p\d)\cap \supp(M)=\emptyset.$ Now using   Proposition \ref{pdelta} for $\fl$-module $M$,  we get $\aa\neq \emptyset.$
\smallskip

\noindent{$\bullet$} $\pmb{\fl\not=A(2k-1,2\ell-1)^{(2)}}:$ Fix $0\neq v\in\aa.$
Suppose that $\dot\a\in \dot R_{ns}^\times.$ Then
there are $\dot\b,\dot\gamma\in \dot R_{sh}$ (see (\ref{new11})) such that $\dot\a=\dot\b+\dot\gamma.$ By Table~\ref{table2},
$$ S_{\dot\b}=S_{\dot\gamma}=\bbbz\d\andd S_{\dot\a}=r\bbbz\d\quad\quad(\hbox{for some $r\in\bbbz^{>0}$}).$$
Since $\boldsymbol\zeta(\d)>0,$ we  choose a large enough $m$ such that $\dot\b+rm'\d,\dot\gamma+r_{\dot \a}m'\d\in R^+\cap R_{re}$ for all $m'\geq m.$
Now as  $v\in \aa,$ for each  nonnegative integer $k$, we have $$\fl^{ \dot\a+r(2m+k)\d}v=[\fl^{ \dot\b+r(m+k)\d},\fl^{ \dot\gamma+rm\d}]v=\{0\}.$$
This completes the proof in this case.

\noindent{$\bullet$} $\pmb{\fl=A(2k-1,2\ell-1)^{(2)}}:$ In this case, $R_{re}\sub R_0.$
Set $$W:=\sum_{\lam\in \supp(M)}\sum_{\ep\in R_{ns}^\times}\fl^\ep M^\lam.$$
Suppose that $\a$ is an element  of the root system $R_0$ of $\fl_0.$ Then $\a$ is either real or imaginary. So if $\ep$ is a nonzero nonsingular root with $\a+\ep\in R,$ we have $\ep+\a\in R_{ns};$ see (\ref{new10}). Therefore,
\begin{align*}
\fl^\a W&=\fl^\a\hbox{\small$\displaystyle{\sum_{\lam\in \hbox{\tiny{\rm supp}}(M)}\sum_{\ep\in R_{ns}^\times}}$}\fl^\ep M^\lam\\&=\hbox{\small$\displaystyle{\sum_{\lam\in \hbox{\tiny{\rm supp}}(M)}\sum_{\ep\in R_{ns}^\times}}$}\fl^\a\fl^\ep M^\lam\\
&\sub \hbox{\small$\displaystyle{\sum_{\lam\in \hbox{\tiny{\rm supp}}(M)}\sum_{\ep\in R_{ns}^\times}}$}\underbrace{[\fl^\a,\fl^\ep]}_{\in \sum_{\eta\in R_{ns}^\times}\fl^\eta} M^\lam+ \hbox{\small$\displaystyle{\sum_{\lam\in \hbox{\tiny{\rm supp}}(M)}\sum_{\ep\in R_{ns}^\times}}$}\fl^\ep\underbrace{\fl^\a M^\lam}_{\in\sum_{\mu\in \hbox{\tiny{\rm supp}}(M)}M^\mu}\\
&\sub W;
\end{align*}
in other words, $W$ is an $\fl_0$-module.
Using Lemma~\ref{new5}, one  finds a positive integer $p$ and $\lam\in \supp(W)$ such that $(\lam+{\bbbz^{>0}}p\d)\cap \supp(W)=\emptyset.$ So by   Proposition \ref{pdelta},
\begin{equation}\label{new9}
\parbox{4.3in}{\small there is a weight $\mu $ of $W$ such that
$\mu+\a$ is not a weight of $W$ if $\a\in  R_0\cap R^+.$}
\end{equation}
 Since $\mu$ is a weight for $W,$  there is a nonzero nonsingular root $\ep$ and $\lam\in \supp(M)$ such that $\fl^\ep M^\lam\neq \{0\}$ and $\mu=\ep+\lam.$ For $0\neq v\in \fl^\ep M^\lam,$ we have
\begin{equation}\label{not-root-w}\fl^\a v\in W^{\a+\mu}\stackrel{(\ref{new9})}{=\joinrel=}\{0\}\quad\quad (\a\in R^+\cap R_0=R^+\cap (R\setminus R_{ns}));
\end{equation} i.e., $v\in\aa.$
We claim that $v$ satisfies (\ref{what-we-need}). We first note that $\dim(\fl^\ep)=1$ and that two times of a {nonzero} nonsingular root is not a root, so
\begin{equation}
\label{g-ep}
\fl^\ep v\in\fl^\ep\fl^\ep M^\lam\sub\underbrace{[\fl^\ep,\fl^\ep]}_{\in\fl^{2\ep}=\{0\}}M^\lam=\{0\}.
\end{equation}
Suppose  $$\ep=\dot\ep+m\d\quad \hbox{for some $\dot\ep\in \dot R^\times_{ns} $ and $m\in \bbbz.$}$$ For each $\dot\a\in \dot R_{ns}^\times,$
by Remark \ref{decom-aff}, one of the following happens:
\begin{itemize}
\item there is $\dot\b_1\in \dot R_{sh}$ such that $\dot \a=\dot\ep+\dot \b_1,$
    \item there are   $\dot\b_1\in \dot R_{sh}$ and $\dot\b_2\in \dot R_{re}^\times$ such that
     $\dot\ep+\dot \b_1\in \dot R_{ns}^\times$ and
     $\dot \a=\dot\ep+\dot \b_1+\dot\b_2,$
     \item there are   $\dot\b_1\in \dot R_{sh}$ and $\dot \b_2,\dot\b_3\in \dot R_{re}^\times$ such that
     $\dot\ep+\dot \b_1,\dot\ep+\dot\b_1+\dot\b_2\in \dot R_{ns}^\times$ and
     $\dot \a=\dot\ep+\dot \b_1+\dot\b_2+\dot\b_3.$
\end{itemize}

$\bullet$ In the first  case, by choosing $t_1\in \bbbz$ with $\boldsymbol{\zeta}(\dot\b_1+t_1\d)>0,$ we have
\begin{align*}
\fl^{\dot \a+t\d}v=[\fl^{\dot\b_1+(t-m)\d},\fl^{ \ep}]v&\sub \fl^{\dot\b_1+(t-m)\d}\fl^{ \ep}v+\fl^{ \ep}\fl^{\dot\b_1+(t-m)\d}v\\
&\stackrel{(\ref{not-root-w}),(\ref{g-ep})}{=\joinrel=\joinrel=\joinrel=}\{0\}\quad\quad (t>t_1+m).
\end{align*}

$\bullet$ In the second case, we
choose $t_1,t_2\in\bbbz^{>0}$ with $t_1+t_2+m>0$ and  $\dot\b_1+t_1\d,\dot\b_2+t_2\d\in R^+.$ Then for $t\geq  t_1+t_2+m,$ by (\ref{not-root-w}), we have  $\fl^{\dot\b_1+(t-t_2-m)\d}v=\{0\}$ and $\fl^{\dot\b_2+t_2\d}v=\{0\}.$  So (\ref{g-ep}) implies that
$$\fl^{\dot \a+t\d}v=[\fl^{\dot\b_2+t_2\d},[\fl^{\dot\b_1+(t-t_2-m)\d},\fl^{ \ep}]]v=\{0\}.$$

$\bullet$ In the third case, we
choose $t_1,t_2,t_3\in\bbbz^{>0}$ with $t_1+t_2+t_3+m>0$ and  $\dot\b_1+t_1\d,\dot\b_2+t_2\d,\dot\b_3+t_3\d\in R^+.$ Then for $t\geq  t_1+t_2+t_3+m,$ as before, we have
$$\fl^{\dot \a+t\d}v=[\fl^{\dot\b_3+t_3\d},[\fl^{\dot\b_2+t_2\d},[\fl^{\dot\b_1+(t-t_2-t_3-m)\d},\fl^{ \ep}]]]v=\{0\}.$$
This completes the proof. \qed

\smallskip

In the following theorem,  we show that the classification problem of  hybrid irreducible finite weight $\fl$-modules  $M$    is reduced to the classification of
cuspidal modules of finite-dimensional cuspidal Levi sub-superalgebras discussed
in \cite{DMP} (see  \cite[Thm. A]{EF} for certain modules over untwisted affine Lie superalgebras).

\begin{Thm}\label{main-1}
Suppose that $M$ is an  hybrid irreducible finite weight $\fl$-module. Then there is a {nontrivial} triangular decomposition $R=R^+\cup R^\circ\cup R^-$ for $R$ and a  triangular decomposition $R^\circ=R^{\circ,+}\cup R^{\circ,\circ}\cup R^{\circ,-}$ for $R^\circ$ with finite $R^{\circ,\circ}$ as well as  a cuspidal finite weight module $N$  over  $\op_{\a\in R^{\circ,\circ}}\fl^\a$ such that $M\simeq{\rm Ind}_\fl(N).$
\end{Thm}
\pf
Suppose that $R=R^+\cup R^\circ\cup R^-$ is the triangular decomposition introduced in the proof of Theorem \ref{triangular}; we mention that $R^\circ$ is finite. We have seen in this theorem that $M^{\fl^+}=\{v\in M\mid \fl^\a v=\{0\}\;\;(\a\in R^+)\}$ is a nonzero module over $\fl^\circ=\op_{\a\in R^\circ}\fl^\a$.
By Proposition \ref{ind}(ii), $M^{\fl^+}$ is an irreducible finite weight $\fl^\circ$-module   and $M\simeq{\rm Ind}_\fl(M^{\fl^+})$. Since $R^\circ$ is finite, $\fl^\circ$ is finite dimensional and so  \cite[Thm. 6.1]{DMP} implies that  there is a triangular decomposition $R^\circ=R^{\circ,+}\cup R^{\circ,\circ}\cup R^{\circ,-}$  for $R^\circ$ and a cuspidal finite weight  module $N$ over  $\op_{\a\in R^{\circ,\circ}}\fl^\a$ such that $M^{\fl^+}\simeq{\rm Ind}_{\fl^\circ}(N).$ This together with Proposition \ref{ind}(ii) and \cite[Cor. 2.4]{DMP} gives that $M\simeq{\rm Ind}_{\fl}(N)$ and so  we are done.
 \qed
\appendix
\section{Affine Lie superalgebras}\label{app}In this section, we recall twisted affine Lie superalgebras from \cite{van-de}.
Suppose that $\fg$ is a finite dimensional  basic classical simple Lie superalgebra with a Cartan subalgebra $\fh\sub\fg_0.$ Suppose that $\kappa$ is  a nondegenerate supersymmetric invariant  even  bilinear form  and $\sg$ is an automorphism of order $n.$ Since $\sg$ preserves $\fg_0$ as well as $\fg_1,$ we have
\begin{align*}
\fg_i=\bigoplus_{k=0}^{n-1}{}^{[k]}\fg_i\quad\hbox{where}\quad {}^{[k]}\fg_i=\{x\in\fg_i\mid \sg(x)=\zeta^kx\}\quad\quad(i\in\bbbz_2,\; 0\leq k\leq n-1).
\end{align*}
in which $\zeta$ is the $n$-th primitive root of unity.  Then
\begin{equation}\label{aff}
\widehat\fg:=\widehat\fg_0\op\widehat\fg_1\quad\hbox{where}\quad\widehat\fg_i=\bigoplus_{k=0}^{n-1}({}^{[k]}\fg_i\ot t^k\bbbc[t^{\pm  n}])\quad\quad(i\in\bbbz_2).
\end{equation}
is a subalgebra of the current superalgebra $\fg\ot\bbbc[t^{\pm1}].$ Setting $$\scg:=\bigoplus_{k=0}^{n-1}({}^{[k]}\fg\ot t^k\bbbc[t^{\pm n}])\op\bbbc c\op\bbbc d\andd \mathscr{H}:=(({}^{[0]}\fg\cap\fh)\ot 1)\op\bbbc c\op\bbbc d. $$ Then $\scg$ together with
$$[x\ot t^p+rc+sd,y\ot t^q+r'c+s'd]:=[x,y]\ot t^{p+q}+p\kappa(x,y)\d_{p+q,0}c+sqy\ot t^q-s'px\ot t^p$$ is a Lie superalgebra called an {\it affine Lie superalgebra} and $\mathscr{H}$ is a Cartan subalgebra of $\scg.$ It is called {\it twisted} if $\sg\neq {\rm id}$ and if $\sg={\rm id}$ and  $\fg\not=A(n,n),$ it is called {\it untwisted}\footnote{The definition of $A(n,n)^{(1)}$ is  slightly different.}. The Lie superalgebra $\scg$ is denoted by $X^{(n)}$ where $X$ is the type of $\fg$.

In what follows, we recall the structure of twisted affine Lie superalgebra  of type $X=A(2k-1,2\ell-1)^{(2)}$ ({\tiny$(k,\ell)\neq (1,1)$}), $A(2k,2\ell)^{(4)},$    $A(2k,2\ell-1)^{(2)}$ and $D(k+1,\ell)^{(2)}$ in which $k,\ell$ are positive integers.

For an integer number $i,$ we define $${\rm sgn}(i):=\left\{
\begin{array}{ll}
1& i>0\\
0& i\leq 0.
\end{array}
\right.$$
For an $m\times n$-matrix $A$ and  positive integers $\ell$ and $k$  define $n\times m$-matrices $A^{\diamond_1},$ $A^{\diamond_2},$ $A^{\diamond_3},$ $A^{\diamond_4}$ and $A^{\diamond_5}$ as follow:
{\small \begin{equation}\label{diamond}
\begin{array}{l}
(A^{\diamond_1})_{r,s}:=(-1)^{r+s}\sg_1(r,s)A_{m+1-s,n+1-r}\\
(A^{\diamond_2})_{r,s}:=(-1)^{r+s}\sg_2(r,s)A_{m+1-s,n+1-r}\quad (\hbox{if }n=2\ell+1)\\
(A^{\diamond_3})_{r,s}:=(-1)^{r+s}\sg_3(r,s)A_{m+1-s,n+1-r}\quad(\hbox{if }m=2\ell+1)\\
(A^{\diamond_4})_{r,s}:=(-1)^{r+s}\sg_4(r,s)A_{m+1-s,n+1-r}\quad (\hbox{if }m=n=2\ell+1)\\
(A^{\diamond_5})_{r,s}:=(-1)^{r+s}\sg_5(r,s)A_{m+1-s,n+1-r}\quad (\hbox{if }n=2k)\\
(A^{\diamond_6})_{r,s}:=(-1)^{r+s}\sg_6(r,s)A_{m+1-s,n+1-r}\quad (\hbox{if }m=2k)\\
(A^{\diamond_7})_{r,s}:=(-1)^{r+s}\sg_7(r,s)A_{m+1-s,n+1-r}\quad (\hbox{if }m=n=2k)
\end{array}
\end{equation}}
where
{\small \begin{equation}\label{sg}
\left\{\begin{array}{rl}
\sg_1(r,s):=&1 \\
\sg_2(r,s):=&(-1)^{{\rm sgn}(r-(\ell+1))}(-1)^{(\ell+1)\d_{r,\ell+1}}i^{\d_{r,\ell+1}}\\
\sg_3(r,s):=&(-1)^{{\rm sgn}(s-(\ell+1))}(-1)^{(\ell+1)\d_{s,\ell+1}}(-i)^{\d_{s,\ell+1}}\\
\sg_4(r,s):=&(-1)^{{\rm sgn}(s-(\ell+1))+{\rm sgn}(r-(\ell+1))}(-1)^{(\ell+1)
(\d_{r,\ell+1}+\d_{s,\ell+1})}i^{\d_{r,\ell+1}}(-i)^{\d_{s,\ell+1}}\\
\sg_5(r,s):=&{(-1)^{{\rm sgn}(k+1-r)}}\\
\sg_6(r,s):=&{(-1)^{{\rm sgn}(k+1-s)}} \\
\sg_7(r,s):=&(-1)^{{\rm sgn}(k+1-r)+{\rm sgn}(k+1-s)}.
\end{array}\right.
\end{equation}}
We note that if $m=n,$ then
\begin{equation}\label{trace}
{\rm tr}(A^{\diamond_1})={\rm tr}(A^{\diamond_4})={\rm tr}(A).
\end{equation}
Also $\diamond_1$ is of order 2 while $\diamond_4$ is of order 4.  Set  $$\fg:=A(m,n)=\left \{\begin{array}{ll}
\frak{psl}(m+1,n+1) & m=n\\
\frak{sl}(m+1,n+1)& m\neq n.
\end{array}
\right.$$
We define
$$h_i:=e_{i,i}-e_{i+1,i+1}\quad d_j:=e_{m+1+j,m+1+j}-e_{m+2+j,m+2+j}\quad (1\leq i\leq m,\;1\leq j\leq n).$$
For $1\leq j\leq m+1$ and $1\leq r\leq n+1,$ define the following functionals on $\fh:=\sspan\{h_i,d_j\mid 1\leq i\leq m,\; 1\leq j\leq n\}$ by
\begin{equation}
\dot\ep_j:\left\{\begin{array}{l}
h_i\mapsto \d_{i,j}-\d_{i+1,j}\\
d_t\mapsto 0
\end{array}
\right.
\quad \dot\d_r:\left\{\begin{array}{l}
h_i\mapsto 0\\
d_t\mapsto \d_{t,r}-\d_{t+1,r}
\end{array}
\right.
\end{equation} for $1\leq i\leq m$ and $1\leq t\leq n.$ The even part  $\fg_0$ of $\fg$ is a reductive Lie algebra  which is centerless if $m=n$ and has a 1-dimensional center  if $m\neq n.$
More precisely, assume
{\small \begin{align}
\ft_1:=&\left\{\left(\begin{array}{cc}
A&{\bf 0}\\
{\bf 0}&{\bf 0}
\end{array}\right)\mid {\rm tr}(A)=0\right\}\simeq\frak{sl}(m+1),
 \ft_2:=\left\{\left(\begin{array}{cc}
{\bf 0}&{\bf 0}\\
{\bf 0}&B\end{array}\right)\mid {\rm tr}(B)=0\right\}\simeq\frak{sl}(n+1),\nonumber\\
\i:=&\left(
\begin{array}{cc}
\frac{1}{m+1}I_{m+1}&{\bf 0}\\
{\bf 0}& \frac{1}{n+1}I_{n+1}
\end{array}
\right).\label{i}
\end{align}}
\hspace{-.1cm}Then the  subalgebras $\fh_1:=\sspan\{h_i\mid 1\leq i\leq m\}$ and $\fh_2:=\sspan\{d_j\mid 1\leq j\leq n\}$ are Cartan subalgebras of $\ft_1$ and $\ft_2$ respectively.
We have  $$\fg_0=\left\{
\begin{array}{lc}
\ft_1\op\ft_2& m=n\\
\ft_1\op\ft_2\op\bbbc\i& m\neq n.
\end{array}
\right.$$
\subsection{$\pmb{A(2k,2\ell)^{(4)}}$} Suppose $m=2k$ and $n=2\ell.$
  For   $X=\left(\begin{array}{cc}
A&B\\
C&D
\end{array}\right)\in \fg,$ define $X^\sg:=\left(\begin{array}{cc}
-A^{\diamond_1}&C^{\diamond_3}\\
-B^{\diamond_2}&-D^{\diamond_4}
\end{array}\right).$ Then $\sg$ defines an automorphism of order 4 on $\fg=A(2k,2\ell).$
The automorphism  $\sg$ maps  each simple component of $\fg_0$ to itself. Suppose $\scg_0(1)$ and $ \scg_0(2)$ are affine Lie algebras obtained from $\ft_1$ and $\ft_2$ using the automorphisms $\sg|_{\ft_1}$ and $\sg|_{\ft_2}$ respectively. Setting
$$\mathscr{H}_i=(({}^{[0]}\fg\cap \fh_i)\ot 1)\op\bbbc c\op\bbbc d\quad(i=1,2),$$ the subalgebra
\begin{equation}\label{cartan1}
\mathscr{H}=
\begin{array}{ll}
\mathscr{H}_1+\mathscr{H}_2
\end{array}
\end{equation}
is a  Cartan subalgebra of $\scg=A(2k,2\ell)^{(4)}$ referred to as the standard Cartan subalgebra.
Contemplating (\ref{i}),  we have $$\scg_0=\left\{\begin{array}{ll}
\scg_0(1)+\scg_0(2)& k=\ell\\
(\scg_0(1)+\scg_0(2))\op(\i\ot t^2 \bbbc[t^{\pm4}])& k\neq \ell.
\end{array}\right.$$   We also  have
$$\scg_0(1)=(\frak{t}_1(\diamond_1)\ot \bbbc[t^{\pm4}])\op (\v\ot t^2\bbbc[t^{\pm4}])\op \bbbc c\op\bbbc d$$ where $\frak{t}_1(\diamond_1)$ and $\v$ are  eigenspaces of  $\ft_1$ corresponding to $1$ and $-1$ respectively with respect to $\diamond_1.$
The automorphism $\diamond_1$ of $\ft_1$ induces an automorphism of the dual space of  $\fh_1,$  mapping $\dot\ep_i-\dot\ep_j$ to $\dot\ep_{2k+2-j}-\dot\ep_{2k+2-i}.$ Setting $\ep_i:=\frac{1}{2}(\dot\ep_i-\dot\ep_{2k+2-i}),$ we get that the set of roots of $\scg_0(1)$ is
$$\mathfrak{R}_1:=(\{\pm\ep_i,\pm\ep_i\pm\ep_j\mid 1\leq i\neq j\leq k\}+2\bbbz\d)\cup (\{\pm2\ep_i\}+4\bbbz\d+2\d)\cup 2\bbbz\d$$ where $\d$ is a functional mapping {$d$ to 1 and $(({}^{[0]}\fg\cap \fh_i)\ot 1)\op\bbbc c$ to $0.$}
Also $\scg_0(2)$ is the affine Lie algebra obtained from {$\frak{t}_2$} by applying $\diamond_4.$ In fact
$$\scg_0(2)=(\frak{t}_2(\diamond_4)\ot \bbbc[t^{\pm4}])\op (\v_\pm\ot t^{\pm 1}\bbbc[t^{\pm4}])\op (\u\ot t^{2}\bbbc[t^{\pm4}]) \op \bbbc c\op\bbbc d$$ where $\frak{t}_2(\diamond_4),$  $\v_\pm$ and $\u$ are  eigenspaces of  $\ft_2$ corresponding to {$1,$} $\pm i$ and $-1$ respectively with respect to $\diamond_4.$
The automorphism $\diamond_4$ induces an automorphism on the dual space of  $\fh_2,$  mapping $\dot\d_j-\dot\d_s$ to $\dot\d_{2\ell+2-s}-\dot\d_{2\ell+2-j}.$ Setting $\d_j:=\frac{1}{2}(\dot\d_j-\dot\d_{2\ell+2-j}),$ we get that the set of roots of $\scg_0(2)$ is
{\small
\begin{align*}
\mathfrak{R}_2&:=(\{\pm2\d_j\mid 1\leq j\leq \ell\}+4\bbbz\d)\cup (\{\pm\d_j\pm\d_s\mid 1\leq j\neq s\leq \ell\}+2\bbbz\d)\\
&\cup (\{\pm\d_i\mid 1\leq i\leq \ell\}+4\bbbz\d\pm\d)\cup 2\bbbz\d.
\end{align*}}

\subsection{$\pmb{A(2k-1,2\ell-1)^{(2)},\;(k,\ell)\neq (1,1)}$} Suppose $m=2k-1$ and $n=2\ell-1.$
For   $X=\left(\begin{array}{cc}
A&B\\
C&D
\end{array}\right)\in \fg,$ define $X^\sg:=\left(\begin{array}{cc}
-A^{\diamond_7}&C^{\diamond_5}\\
-B^{\diamond_6}&-D^{\diamond_1}
\end{array}\right).$ Then $\sg$ defines an automorphism of order 2 on $\fg=A(2k-1,2\ell-1).$ Set  $\scg=A(2k-1,2\ell-1)^{(2)}$ and suppose $\scg_0(1)$ and $\scg_0(2)$ are affine Lie algebras obtained by {the affinization of} $\ft_1$ and $\ft_2$ using the automorphism $\sg.$  Then we have
$$\scg_0=\left\{
\begin{array}{ll}
\scg_0(1)+\scg_0(2)& k=\ell\\
(\scg_0(1)+\scg_0(2))\op(\i\ot t\bbbc[t^{\pm2}])&k\neq \ell;
\end{array}
\right.
$$ see (\ref{i}). Setting
$$\mathscr{H}_i=(({}^{[0]}\fg\cap \fh_i)\ot 1)\op\bbbc c\op\bbbc d\quad(i=1,2),$$ we get that
\begin{equation}\label{cartan1}
\mathscr{H}=
\begin{array}{ll}
\mathscr{H}_1+\mathscr{H}_2
\end{array}
\end{equation}
is a  Cartan subalgebra of $\scg.$ We call it  the standard Cartan subalgebra of $\scg$.
 We have
$$\scg_0(1)=(\frak{t}_1(\diamond_7)\ot \bbbc[t^{\pm2}])\op (\v\ot t\bbbc[t^{\pm2}])\op \bbbc c\op\bbbc d$$ where $\frak{t}_1(\diamond_7)$ and $\v$ are  eigenspaces of  $\ft_1$ corresponding to $1$ and $-1$ respectively with respect to $\diamond_7.$
The automorphism $\diamond_7$ of $\ft_1$ induces an automorphism of the dual space of  $\fh_1,$  mapping $\dot\ep_i-\dot\ep_j$ to $\dot\ep_{2k+1-j}-\dot\ep_{2k+1-i}.$ Setting $\ep_i:=\frac{1}{2}(\dot\ep_i-\dot\ep_{2k+1-i}),$ we get that the set of roots of $\scg_0(1)$ is
$$\mathfrak{R}_1:=(\{\pm\ep_i\pm\ep_j\mid 1\leq i\neq j\leq k\}+\bbbz\d)\cup (\{\pm2\ep_i\}+2\bbbz\d+\d)\cup \bbbz\d$$ where $\d$ is a functional mapping {$d$ to 1 and $(({}^{[0]}\fg\cap \fh_i)\ot 1)\op\bbbc c$ to $0.$}
Also $\scg_0(2)$ is the affine Lie algebra obtained from $\frak{t}$ by applying $\diamond_1.$ In fact
$$\scg_0(2)=(\frak{t}_2(\diamond_1)\ot \bbbc[t^{\pm2}])\op (\v\ot {t}\bbbc[t^{\pm2}]) \op \bbbc c\op\bbbc d$$ where $\frak{t}_2(\diamond_1)$ and   $\v$ are  eigenspaces of  $\ft_2$ corresponding to $1$ and $-1$ respectively with respect to $\diamond_1.$
The automorphism $\diamond_1$ induces an automorphism on the dual space of  $\fh_2$ mapping $\dot\d_j-\dot\d_s$ to $\dot\d_{2\ell+1-s}-\dot\d_{2\ell+1-j}.$ Setting $\d_j:=\frac{1}{2}(\dot\d_j-\dot\d_{2\ell+1-j}),$ we get that the set of roots of $\scg_0(2)$ is
{\small$$\mathfrak{R}_2:=(\{\pm2\d_j\mid 1\leq j\leq \ell\}+2\bbbz\d)\cup (\{\pm\d_j\pm\d_s\mid 1\leq j\neq s\leq\ell\}+\bbbz\d)\cup \bbbz\d.$$}

\subsection{$\pmb{A(2k,2\ell-1)^{(2)}}$} Suppose $m=2k$ and $n=2\ell-1.$
For   $X=\left(\begin{array}{cc}
A&B\\
C&D
\end{array}\right)\in \fg,$ define $X^\sg:=\left(\begin{array}{cc}
-A^{\diamond_1}&C^{\diamond_1}\\
-B^{\diamond_1}&-D^{\diamond_1}
\end{array}\right).$ Then $\sg$ defines an automorphism of order 2 on $\fg=A(2k,2\ell-1).$
For $\scg=A(2k,2\ell-1)^{(2)},$ the Cartan subalgebra of $\scg$ is
\begin{equation}\label{cartan1}
\mathscr{H}=\mathscr{H}_1+\mathscr{H}_2 \hbox{ with $\mathscr{H}_i=(({}^{[0]}\fg\cap \fh_i)\ot 1)\op\bbbc c\op\bbbc d\quad(i=1,2).$}
\end{equation}
The  Cartan subalgebra $\hH$ is called the standard Cartan subalgebra of $\scg.$
Moreover,  we have $\scg_0=\scg_0(1)+\scg_0(2)\op(\i\ot t\bbbc[t^{\pm2}]),$ where $\scg_0(1)$ is the affine Lie algebra obtained from $\frak{t}_1$ by applying $\diamond_1;$ in fact
$$\scg_0(1)=(\frak{t}_1(\diamond_1)\ot \bbbc[t^{\pm2}])\op (\v\ot t\bbbc[t^{\pm2}])\op \bbbc c\op\bbbc d$$ where $\frak{t}_1(\diamond_1)$ and $\v$ are  eigenspaces of  $\ft_1$ corresponding to $1$ and $-1$ respectively with respect to $\diamond_1.$
The automorphism $\diamond_1$ of $\ft_1$ induces an automorphism of the dual space of  $\fh_1$  mapping $\dot\ep_i-\dot\ep_j$ to $\dot\ep_{2k+2-j}-\dot\ep_{2k+2-i}.$ Setting $\ep_i:=\frac{1}{2}(\dot\ep_i-\dot\ep_{2k+2-i}),$ we get that the set of roots of $\scg_0(1)$ is
$$\mathfrak{R}_1:=(\{\pm\ep_i,\pm\ep_i\pm\ep_j\mid 1\leq i\neq j\leq k\}+\bbbz\d)\cup (\{\pm2\ep_i\}+2\bbbz\d+\d)\cup \bbbz\d$$ where $\d$ is a functional mapping {$d$ to 1 and $(({}^{[0]}\fg\cap \fh_i)\ot 1)\op\bbbc c$ to $0.$}
Also $\scg_0(2)$ is the affine Lie algebra obtained from $\frak{t}$ by applying $\diamond_1.$ In fact
$$\scg_0(2)=(\frak{t}_2(\diamond_1)\ot \bbbc[t^{\pm2}])\op (\v\ot {t}\bbbc[t^{\pm2}]) \op \bbbc c\op\bbbc d$$ where $\frak{t}_2(\diamond_1)$ and   $\v$ are  eigenspaces of  $\ft_2$ corresponding to $1$ and $-1$ respectively with respect to $\diamond_1.$
The automorphism $\diamond_1$ induces an automorphism on the dual space of  $\fh_2,$ consisting of all diagonal matrices, mapping $\dot\d_j-\dot\d_s$ to $\dot\d_{2\ell+1-s}-\dot\d_{2\ell+1-j}.$ Setting $\d_j:=\frac{1}{2}(\dot\d_j-\dot\d_{2\ell+1-j}),$ we get that the set of roots of $\scg_0(2)$ is
{\small$$\mathfrak{R}_2:=(\{\pm2\d_j\mid 1\leq j\leq \ell\}+2\bbbz\d)\cup (\{\pm\d_j\pm\d_s\mid 1\leq j\neq s\leq \ell\}+\bbbz\d)\cup \bbbz\d.$$}

\subsection{$\pmb{D(k+1,\ell)^{(2)}}$}
 We know that $\fg:=\frak{osp}(2k+2,2\ell)$ consists of all matrices of the form{
\begin{equation}\label{matrix form}
\left(
\begin{array}{rcc}
\begin{array}{cc}
x&y\\
z&-x^t
\end{array}
& \vline&
\begin{array}{cc}
m&n\\
p&q
\end{array}\\
\hline
\begin{array}{cc}
-q^t&-n^t\\
p^t&m^t
\end{array}&\vline& \begin{array}{cc}
r&s\\
u&-r^t
\end{array}
\\
\end{array}\right)
\end{equation}
where $x,$ $m$ and $r$ are respectively  $(k+1)\times (k+1)$, $(k+1)\times \ell$ and $\ell\times\ell$-matrices and $y$ as well as $z$ are skew-symmetric matrices while $s$ and $u$ are symmetric.} We make  a convention that for $1\leq i\leq k+1,$  set $\bar i:=i+k+1.$
Set $G:=(g_{i,j})$ to be a $(2k+2)\times(2k+2)$-matrix define dy $$g_{i,j}=g_{\bar i,\bar j}:=(1-\d_{i,k+1})\d_{i,j}\andd g_{\bar i,j}=g_{j,\bar i}=\d_{i,k+1}\d_{i,j}\quad\quad(1\leq i,j\leq k+1).$$
Then $G$ is invertible with $G^{-1}=G.$ Next set $$H:=\left(
\begin{array}{cc}
  G&{\pmb 0}\\
  {\pmb 0}& I_{2\ell}
\end{array}
\right)$$
in which $I_{2\ell}$ is the identity matrix of dimension $2\ell.$ The automorphism $\sg$ mapping $X\in \fg$ to $HXH^{-1}$ is an automorphism of $\fg$ of order 2. We have $\fg_0=\ft_1\op\ft_2$ where
\begin{align*}
\ft_1\simeq D(k+1)\andd
 \ft_2\simeq C({\ell}).
\end{align*}
In fact $\ft_1$ (resp. $\ft_2$) consists of block matrices of the form (\ref{matrix form}) whose second, third and fourth (resp. first) block are zero  matrices.
Suppose $\fh_1$ is the abelien subalgebra of $\ft_1$ spanned by $\{h_i:=e_{i,i}-e_{\bar i,\bar i}\mid 1\leq i\leq k+1\}$ and  $\fh_2$ is the abelien subalgebra of $\ft_2$ spanned by $\{d_p:=e_{2k+2+p,2k+2+p}-e_{2k+2+\ell+p,2k+2+\ell+p}\mid 1\leq p\leq \ell\}.$
Define
\begin{align*}
 \ep_i:&\fh_1^*\longrightarrow \bbbc& \d_p:&\fh_2^*\longrightarrow \bbbc\\
 &h_j\mapsto \d_{i,j}&&d_q\mapsto \d_{p,q}
 \end{align*}
where $1\leq i,j\leq k+1$ and $1\leq p,q\leq \ell$. Then $\{\ep_i\mid 1\leq i\leq k+1\}$ is a basis for the dual space $\fh_1^*$ of $\fh_1$ and $\{\d_p\mid 1\leq p\leq \ell\}$ is a basis for the dual space $\fh_2^*$ of $\fh_2.$

For $\scg=D(k+1,\ell)^{(2)},$
the standard Cartan subalgebra of $\scg$ is
\begin{equation}\label{cartan1}
\mathscr{H}=\mathscr{H}_1+\mathscr{H}_2 \hbox{ with $\mathscr{H}_i=(({}^{[0]}\fg\cap \fh_i)\ot 1)\op\bbbc c\op\bbbc d\quad(i=1,2).$}
\end{equation}
Moreover,  we have  $\scg_0=\scg_0(1)+\scg_0(2)$ where $\scg_0(i)$ ($i=1,2$) is the affine Lie algebra obtained from $\frak{t}_i$ by applying $\sg.$ In fact
$$\scg_0(1)=(\frak{t}_1(\sg)\ot \bbbc[t^{\pm2}])\op (\v\ot t\bbbc[t^{\pm2}])\op \bbbc c\op\bbbc d$$ where $\frak{t}_1(\sg)$ and $\v$ are  eigenspaces of  $\ft_1$ corresponding to $1$ and $-1$ respectively with respect to $\sg|_{\frak{t}_1}.$

The automorphism $\sg$  induces an automorphism of the dual space of $\fh_1^*$ mapping
 $$\ep_i\mapsto \left\{
\begin{array}{cc}
\ep_i& i\in\{1,\ldots,k\}\\
-\ep_{i}& i=k+1.
\end{array}
\right.$$
 The set of roots of $\scg_0(1)$ is
$$\mathfrak{R}_1:=(\{0,\pm\ep_i,\pm\ep_i\pm\ep_j\mid 1\leq i\neq j\leq k\}+2\bbbz\d)\cup (\{0,\pm\ep_i\}+2\bbbz\d+\d)$$  where $\d$ is a functional mapping {$d$ to 1 and $(({}^{[0]}\fg\cap \fh_i)\ot 1)\op\bbbc c$ to $0.$}
The automorphism $\sg$ is the identity map on $\ft_2$ and so
$$\scg_0(2)=(\frak{t}_2\ot \bbbc[t^{\pm2}]) \op \bbbc c\op\bbbc d.$$ The root system of $\scg_0(2)$ is
$$\mathfrak{R}_2=\{\pm\d_p\pm\d_q\mid 1\leq p,q\leq \ell\}+2\bbbz\d.$$

\begin{rem}\label{decom-aff-app}
{\rm
 As we have seen if $\scg=\scg_0\op\scg_1$ is a twisted affine Lie superalgebra  of type $X=A(2k-1,2\ell-1)^{(2)}$ ({\tiny$(k,\ell)\neq (1,1)$}), $A(2k,2\ell)^{(4)},$    $A(2k,2\ell-1)^{(2)}$ and $D(k+1,\ell)^{(2)}$ where  $k,\ell$ are positive integers,
there are  affine Lie subalgebras $\scg_0(1)$ and $\scg_0(2)$ of $\scg_0$ with Cartan subalgebras $\mathscr{H}_1$ and $\mathscr{H}_2$ respectively such that $$\mathscr{H}:=\mathscr{H}_1+\mathscr{H}_2$$ is a Cartan subalgebra of $\scg$ and up to an $\hH$-module whose weights are  nonzero imaginary roots, $\scg_0$ equals  {$\scg_0(1)\op\scg_0(2)$.}
}
\end{rem}
\centerline{\bf Acknowledgment}
The  author would like to express her thanks to Alberto Elduque for some helpful discussion and other people of Department of Mathematics, University of Zaragoza, for their kind hospitality during her visit when some part of this work has been done. She also would like to thank
Michael Lau for his comments and introducing  \cite{DG} to her
and Karl-Hermann Neeb
for his comments on this paper. At the end she thanks the anonymous referee for his/her comments and suggestions on the paper.

\end{document}